\documentclass[10pt]{article}
\usepackage{authblk}
\usepackage{latexsym}
\usepackage{amssymb,amsbsy,amsmath,amsfonts,amssymb,amscd}
\usepackage{epsfig, graphicx}
\usepackage{mathtools}
\usepackage{mathrsfs}
\numberwithin{equation}{section}
\usepackage{colortbl}
\newcommand{\commentout}[1]{}
\newcommand{\C}{\mathbb{C}}
\newcommand{\R}{\mathbb{R}}
\newcommand{\N}{\mathbb{N}}

\newcommand {\Chi} {{\bf \raise 2pt \hbox{$\chi$}} }
\newcommand {\f}   {\frac}
\newcommand {\p}   {\partial}

\newcommand {\proof} {\noindent {\bf Proof}. }
\newcommand{\bqn}{\begin{eqnarray}}
\newcommand{\eqn}{\end{eqnarray}}
\newcommand{\bqns}{\begin{eqnarray*}}
\newcommand{\eqns}{\end{eqnarray*}}
\newcommand{\beq}{\begin{equation}}
\newcommand{\beqa} {\begin{array}{rl}}
\newcommand{\eeq}{\end{equation}}
\newcommand{\eeqa}{\end{array}}
\newtheorem{theorem}{Theorem}[section]

\newtheorem{remark}[theorem]{Remark}
\newtheorem{proposition}[theorem]{Proposition}
\newtheorem{corollary}[theorem]{Corollary}
\newcommand{\qed}{{ \hfill
                       {\unskip\kern 6pt\penalty 500
                       \raise -2pt\hbox{\vrule\vbox to 6pt{\hrule width 6pt
                       \vfill\hrule}\vrule} \par}   }}

\title{On the non existence of  non negative solutions to a critical Growth-Fragmentation Equation}

\author{
 M.  Escobedo \thanks{{Departamento de Matem\'aticas, Universidad del Pa{\'\i}s Vasco (UPV/EHU). E--48080 Bilbao (SPAIN).\hfill \break  email: miguel.escobedo@ehu.es}}
}
\date{}

\begin{document}
\maketitle
\pagestyle{plain}
%\tableofcontents
\pagenumbering{arabic}

\begin{abstract}
A  growth fragmentation equation with constant dislocation density measure is considered, in which  growth and division rates balance each other. This  leads to a simple  example of equation where the so called Malthusian hypothesis $(M_+)$ of 
J. Bertoin and A. Watson (2016) is not necessarily satisfied. It is proved that when  that happens, and as it was first suggested by these authors, no global non negative weak solution, satisfying some boundedness condition on several of its moments, exist. Non existence of  local non negative solutions satisfying a similar condition, is proved to happen also. When a local non negative solution exists, the explicit expression is given.  
\end{abstract}

%\tableofcontents
%%%%%%%%%%%%%%%%%%%%%%%%%%%%%%%%%%%%%%%%%%%
\section{Introduction}
Growth fragmentation equations  have  proved to be  of interest due to  their many  applications in mathematical modeling  and also for purely mathematical reasons (cf. \cite{MMP}, \cite{BA}, \cite{DG}, \cite{BG} and references therein).  Motivated by the study of compensated growth-fragmentation stochastic processes (cf. \cite{B}) and their occurrence in the construction of the Brownian map (cf. \cite{BBCK}, \cite{LG}, \cite {M}), the  Cauchy problem for the equation 
\bqn
\label{eq:croisfrag2}
&&\f{\p}{\p t} u(t,x) + \f{\p}{\p x} \big(x^{1+\gamma } u (t,x)\big) +  x^\gamma u(t,x)=\int_x^\infty \frac {1} {y} k_0\left(\frac {x} {y} \right)y^\gamma   u(t,y) dy,\,\,\,t>0, x>0
\eqn
is considered in \cite{BW} with  initial data:
\bqn
u(0)=\delta_1\label{eq:croisfrag2data}, 
\eqn
for $\gamma\in \R $ and  $k_0$ a  dislocation measure density, with support contained in $[0, 1] $ and satisfying: 
\bqn
\label{S0E1}
k_0(x)dx=k_0(1-x)dx,\,\,\,\forall x\in [1/2, 1); \quad\int  _{ [1/2, 1) }\!\!\!\!(1-x)^2k_0(x)dx<\infty.
\eqn
The existence of solutions of growth fragmentation equations has been studied by several authors, with different motivations, by different methods, for different dislocation measures, and fragmentation rates,  (cf. for example \cite{MMP}, \cite{BA}, \cite{DG}, \cite{BG} and references therein). However, the equation  \eqref{eq:croisfrag2} is rather specific. It is said critical because the growth and the dislocation rates balance each other. The case  $\gamma =0$ was considered in  \cite{DE}, although for less general dislocation measures. When  $\gamma \not =0$  the growth rate is not linear  and  the dislocation rate unbounded or singular. In that case, the existence of  global, non negative, weak solutions of \eqref{eq:croisfrag2}-\eqref{S0E1} has been  proved in \cite{BW}, under the  condition (called Malthusian condition $(M_+)$ in \cite{BW}):
\bqn
&& \inf _{ s\ge 0 }\Phi (s)<0 \label{S1prop1B}\\
&&\Phi (s)= \left(K(s)+s-2\right),\,\,\,K(s)=\int _0^\infty x^{s-1}k_0(s)ds.\label{S1prop2}
\eqn

When property \eqref{S1prop1B} is not satisfied it is  shown in \cite{BW} and \cite{BS} that  the particle system that corresponds to the stochastic version of 
\eqref{eq:croisfrag2}-\eqref{S0E1} explode in finite time almost surely.
The question has then been raised in  \cite{BW} of the existence of  non negative global  solutions to  \eqref{eq:croisfrag2}-\eqref{S0E1} when the measure $k_0$ is such that:
\bqn
\inf _{ s\ge 0 } \Phi (s)\ge 0 \label{S1prop1}
\eqn
and it was suggested that no such solutions exists when the inequality in \eqref{S1prop1} is strict. In order to obtain some insight into this question, we consider the simplest possible choice for $k_0$:
\bqn
\label{S1EK0}
k_0(x)=\theta H(1-x),\,\,\,\, \theta>0
\eqn
where $H$ is the Heaviside's function.  This is of course a very particular example, but  for which it is possible to obtain a rather explicit solutions, whose properties may be understood in detail. It is straightforward to check that for such a dislocation measure:
\bqn
&&\mathcal M _{ k_0 }(s)=\frac {\theta} {s}\,\,\,\hbox{and}\,\,\,\,\Phi (s)=\frac {\theta} {s}+s-2\equiv \frac{(s-\sigma _1)(s-\sigma _2)}{s},\,\,\forall s\in \C;\,\,\Re e(s)>0, \label{S1EKS0}\\
&&\sigma _1=1-\sqrt{1-\theta},\,\,\,\sigma _2=1+\sqrt{1-\theta}.\label{S1EKS024}
\eqn
If $\theta\in (0, 1)$, the two roots of $\Phi (s)$ are positive real numbers and  condition  (\ref{S1prop1B}) is satisfied. 
But, when $\theta \ge 1$,  $\inf _{ s>0 } \Phi (s)=2(\sqrt \theta -1)\ge 0$ and  (\ref{S1prop1B}) is not satisfied.  

For $\theta\in (0, 1)$ the existence of global non negative  solutions  follows from the results of \cite{BW}. We then focus on the case $\gamma \not=0$, $\theta>1$ and the question of the existence or not of  non negative solutions.
\subsection{Some notations.}
We denote $\N$ the set of non negative integers and $\Gamma (\cdot)$ the Gamma function. For  a given interval $(a, b)\subset \R$ we define:
\bqn
\label{S4ES}
\mathscr{S}(a, b)=\left\{s\in \C;\,\,\,\Re e(s)\in (a, b) \right\}.
\eqn
We denote $\mathscr{D}_1'$  the set of distributions of order one and by $F(a, b, c, z)$  the Gauss hypergeometric function $_2F_1(a, b; c; z)$.
We  say that the measure $u$ is a weak solution of \eqref{eq:croisfrag2},\eqref{S1EK0} on the time interval $(t_0, t_1)$ if
\bqn
&&\forall \varphi \in C_c^1((t_0, t_1)\times (0, \infty)):\nonumber\\
&&\int _0^\infty\int _0^\infty\left(\frac {\partial \varphi } {\partial t}+x^{\gamma +1}\frac {\partial \varphi } {\partial  x}+x^\gamma \varphi  (t, x) \right)u(t, x)dxdt=\nonumber\\
&&\hskip 5cm -\theta \int _0^\infty \int _0^\infty u(t, y)y^{\gamma -1}\int _0^y \varphi (t, x)dxdydt \label{S1Esol21}
\eqn
We denote $\mathscr{M} _{ \rho  }$ the space of measures $u$ on $(0, \infty)$ such that
\bqns
\int _0^\infty x^{\rho } u(x)dx<\infty.
\eqns
If  $w$ is a measure, we denote $\mathcal M _{ w }$ its Mellin transform, defined, when it makes sense, as
$$
\mathcal M_u(t, s)=\int _0^\infty x^{s-1}u(t, x)dx.
$$
It follows from the definition of $K(s)$ in \eqref{S1prop2} that $K(s)=\mathcal M _{ k_0 }(s)$. The use of the Mellin transform makes the spaces $E' _{ p,\,q }$ for $p<q$, presented for example in Chapter 11 of \cite{ML}, necessary. They are  defined as the dual of the spaces $E _{ p,\,q }$ of all the functions $\phi \in \mathscr{C}^\infty (0, \infty)$ such that:
\bqns
N _{ p, q, k }(\phi )=\sup _{ x>0 }\left(k _{ p, q }(x)x^{k+1}\left|\phi ^{k}(x) \right|\right)<\infty,\,\,\hbox{where}\,\,\,
k _{ p, q }(x)=
\left\{
\begin{split}
&x^{-p},\,\,\hbox{if}\,\,0<x\le 1\\
&x^{-q},\,\,\hbox{if}\,\,x> 1
\end{split}
\right.
\eqns
with the topology defined by the numerable set of seminorms $\left\{N _{ p, q, k } \right\} _{ k\in \N }$. It follows that $E' _{ p,\,q }$ is a subspace of $\mathscr{D}'(0, \infty)$.  As indicated  in \cite{ML}, these are the spaces of Mellin transformable distributions.
\subsection{Main results.}

In very short, when $\theta >1$ and $\gamma \not =0$, global  non negative solutions to \eqref{eq:croisfrag2},\eqref{S1EK0},\eqref{eq:croisfrag2data}, satisfying a  boundedness condition on several of its moments, do not exist. But more detailed statements depend on the sign of $\gamma $, as follows.
\subsubsection{When $\gamma >0$.}
Our first  result is the following local existence of non negative solutions when $\gamma >0$:
\begin{theorem}
\label{S1MainTh1}
For all  $\theta>0$ and $\gamma >0$ there exists a unique, non negative weak solution $u \in \mathscr{D}_1'((0, \gamma ^{-1})\times (0, \infty))$ of  \eqref{eq:croisfrag2},\eqref{S1EK0} on $(0, \gamma ^{-1})$ such that, for some $\rho >0$:
\bqn
&&u \in C\left(\left[0, \gamma ^{-1}\right); E' _{ \rho -\delta , \rho +\gamma+ \delta  } \right),\,\,\,\hbox{for some}\,\,\,\delta >0 \label{S1E40987}\\
%&& \forall T\in \left(0, \gamma ^{-1}\right),\,\,\exists C_T>0;\,\, \int _0^\infty u(t, x)x^{s-1} dx\le C_T\,\,\,\,\forall t\in [0, T], \,\,\forall s\in [\rho , \rho +\gamma ]\label{S1E41789}\\
&&u(t)\rightharpoonup \delta_1\,\,\,\hbox{in}\,\,E' _{ \rho , \rho +\gamma  },\,\,\, \hbox{as}\,\,t\to 0.\label{S1E42}
\eqn
That solution is:
\bqn
u (t, x)&=&u ^S(t, x)+u ^R(t, x)H\left(1-\left(1-\gamma  t\right)^{\frac {1} {\gamma }}x\right) \label{S1Esol0}\\
u ^S(t, x)&=&\left(1-\gamma  t\right)^{\frac {1} {\gamma }}\delta \left(x-\left(1-\gamma  t\right)^{-\frac {1} {\gamma }} \right)\label{S1Esol00}\\
u ^R(t, x)&=&\theta\left(1-\gamma  t\right)^{\frac {2} {\gamma }}tF\left(1+\frac {\sigma _1} {\gamma }, 1+\frac {\sigma _2} {\gamma }, 2, \gamma  t\left(1+\left(\gamma  t-1\right)x^\gamma  \right) \right)
\label{S1Esol}\\
&&\hskip -4.8cm \hbox{and satisfies:}\hskip 3.2cm 
u\in \mathscr{C}\left(\left[0, \gamma ^{-1}\right); \mathscr{M} _{ \rho-1  } \right),\,\,\forall \rho >0. \label{S1Esol1000}
\eqn
\end{theorem}

The sense in which the initial data $\delta_1$ is taken in the hypothesis \eqref{S1E42}  ensures that the Mellin transform of $u(t)$ converges to $1$ as $t$ goes to zero, for all $s\in \mathscr{S}(\rho , \rho+ \gamma )$. Since $E' _{ \rho , \rho +\gamma  }\subset \mathscr{D}'(0, \infty)$ with continuous embedding, this condition is stronger than the convergence in the weak sense of measures.

Non uniqueness  in some sense, of non negative solutions of \eqref{eq:croisfrag2} has been proved in \cite{BW} under some conditions on $k_0$. However the function $k_0$ given  in \eqref{S1EK0}  does not satisfy such conditions (cf. Remark \ref{S3Run}). 

As a consequence of Theorem \ref{S1MainTh1} we deduce the following result:
\begin{corollary}
\label{S1Cor}
The  solution $u$ of  \eqref{eq:croisfrag2},\eqref{S1EK0} defined in \eqref{S1Esol0}-\eqref{S1Esol} satisfies:
\bqn
\label{S1Ebw0}
\lim _{ \gamma  t\to 1^- }u (t, x)=\frac {\gamma \Gamma \left(\frac {2} {\gamma } \right)}
{  \Gamma \left(\frac {\sigma _1} {\gamma }  \right) \Gamma \left(\frac {\sigma _2} {\gamma }  \right)}(1+x^\gamma )^{-\frac {2} {\gamma }},\,\,\,\forall x>0
\eqn
and:
\bqn
\forall r>1:\,\,\, \lim _{ t\to \gamma ^{-1} }(1-\gamma t)^{\frac {r-1} {\gamma }}\int _0^\infty x^r u (t, x)dx&=&\frac {\Gamma \left(\frac {r+1} {\gamma }\right)\Gamma \left(\frac {r-1} {\gamma }\right)}
{\Gamma \left(\frac {r+1-\sigma _1} {\gamma }\right)\Gamma \left(\frac {r+1-\sigma _2} {\gamma }\right)},\label{S1Ebw1} \\
 \lim _{ t\to \gamma ^{-1} }\frac {-1} {\log (1-\gamma t)}\int _0^\infty x u(t, x)dx&=&\frac {\Gamma \left(\frac {2} {\gamma }\right)}
{\Gamma \left(\frac {\sigma _1} {\gamma }\right)\Gamma \left(\frac {\sigma _2} {\gamma }\right)},\label{S1Ebw2}\\
\forall r\in (0, 1):\,\,\,\,
\lim _{ t\to \gamma ^{-1} }\int _0^\infty x^r u(t, x)dx&=&\frac {\Gamma \left(\frac {r+1} {\gamma }\right)\Gamma \left(\frac {1-r} {\gamma }\right)}
{\Gamma \left(\frac {\sigma _1} {\gamma }\right)\Gamma \left(\frac {\sigma _2} {\gamma }\right)}.\label{S1Ebw3}
\eqn
\end{corollary}

The two properties \eqref{S1Ebw1} and \eqref{S1Ebw2} show that the moments of order $r\ge 1$ of  $u$, that are finite for all $\gamma t<1$, become infinite as $\gamma t\to 1^-$. 

We prove in Theorem \ref{S5T621} that, when $\theta>1$ and $\gamma > 0$, there is no possible extension of $u$ to a non negative global solution whose Mellin transform satisfies  suitable conditions. When $\gamma \in (0, 2)$ and $\theta>1$ the non existence of non negative  solutions for large times  is shown in the following:
\begin{theorem}
\label{S1MainTh4Bis} Suppose that $\gamma \in (0, 2)$, $\theta>1$ and $T>\gamma ^{-1}$. Then, there is no  extension of the local solution $u$ to a non negative weak solution $w\in 
 \mathscr{D}_1'((0, T)\times (0, \infty))$ of  \eqref{eq:croisfrag2},\eqref{S1EK0} such that   for some $\rho\in (0, 2-\gamma )$:
\bqn
&&w \in C\left(\left[0, T\right); \mathscr{M} _{ \rho-1 -\delta  }\cap\mathscr{M} _{ \rho+\gamma -1 +\delta  } \right),\,\,\,\hbox{for some}\,\,\,\delta >0 \label{S1E40}
\eqn
and  satisfying the initial condition \eqref{S1E42}.
\end{theorem}

\subsubsection{When $\gamma <0$.}
When $\gamma <0$ the existence of a local solution $v\in \mathscr C([0, -\gamma ^{-1}), E' _{ 1+\gamma , \infty })$ of  \eqref{eq:croisfrag2},\eqref{S1EK0},\eqref{S1E42}  on $(0, -\gamma ^{-1})$  is proved in Theorem \ref{S6T12345}.  But the following non existence of local nonnegative solutions holds:
\begin{theorem}
\label{S1Th421}
If $\gamma <0$ and $\theta >1$ there is no local, non negative weak solution $v$ of  \eqref{eq:croisfrag2},\eqref{S1EK0} on $(0, T)$,  for any $T>0$, satisfying for some $\rho>1-\gamma $,   
the initial condition \eqref{S1E42} and such that: 
\bqn
v  \in C\left(\left[0, T\right); \mathscr{M} _{ \rho-1 -\delta  }\cap\mathscr{M} _{ \rho-\gamma -1 +\delta  } \right),\,\,\,\hbox{for some}\,\,\,\delta >0. \label{S1ENeg40}
\eqn
\end{theorem}

In summary,  for $\theta>1$ and $\gamma \in (0, 2)$, the non negative solution $u$ of Theorem \ref{S1MainTh1} blows up as $\gamma t\to 1^-$, in the sense given by Corollary  \ref{S1Cor}, and can not be extended beyond $t=\gamma ^{-1}$ to a non negative solution  that still satisfies \eqref{S1E40}
. If $\gamma <0$, nonnegative solutions satisfying  \eqref{S1ENeg40}
do not exist, even locally in time. 
 The non existence of global non negative solutions, for  critical growth fragmentation equations where the condition \eqref{S1prop1B} is not satisfied, was first suggested  in~\cite{BW}. Of course, Theorem \ref{S1MainTh4Bis} and Theorem \ref{S1Th421} do not preclude the existence of  non negative global solutions that do not satisfy  \eqref{S1E40} or \eqref{S1ENeg40}
. When $\theta\in (0, 1)$, the condition \eqref{S1prop1B} is satisfied and then, as  proved in \cite{BW}, the  problem \eqref{eq:croisfrag2},\eqref{eq:croisfrag2data}  has a global non negative solution $\mu $. It  follows that $\mu$  coincides with the solution  obtained in Section  \ref{sign}, when  $\gamma<\sqrt{1-\theta}$ (cf. Proposition \ref{finalrem}) or,  when $\gamma <0$, with that obtained in Section \ref{SGammaNegatif} (cf. Remark \ref{S6R1}). If $\theta=1$ the  condition \eqref{S1prop1B}  is not satisfied, but our arguments do not prove the non existence of a non negative extension of $u$ beyond $t=\gamma ^{-1}$ (cf. Remark \ref{S5Etheta1}).

The equation  \eqref{eq:croisfrag2} may be solved taking advantage of its linearity, using the  Mellin  transform. The proof of the non existence of  non negative solution is then done in two steps. The first is to prove the uniqueness of solutions that may take positive and negative values but some moments of which are suitably bounded. The second is to show that the solution that was previously obtained satisfies the regularity condition, but takes  positive and negative values. That follows from its behavior as $x\to 0$ or $x\to \infty$, since it is  given, up to some multiplicative factor depending on time,  by $x^{-\sigma _2-\gamma }$ and  $x^{-\sigma _1-\gamma }$. When $\sigma _2$ and $\sigma _1$ are complex numbers, this forces the solution to oscillate.

The  choice of $k_0$ as in \eqref{S1EK0}, is of course very particular and makes  the solutions of equation \eqref{eq:croisfrag2} rather explicit. We may recall  at this point  that explicit solutions to the  Cauchy problem for the  pure fragmentation equation (i.e. without   growth term and with  $k_0$ such that  $\int y k_0(y)dy=1$), with  the initial data as in \eqref{eq:croisfrag2data}  where obtained in \cite{ZG1}, \cite{ZG2} for several fragmentation rates and  the same dislocation measure \eqref{S1EK0} with $\theta=2$. We emphasize however that the arguments used in Section~\ref{Extension} and Section \ref{SGammaNegatif}, based on the Wiener Hopf method,  permit to solve  the growth fragmentation equation \eqref{eq:croisfrag2} for more general dislocation measures. More details will be presented elsewhere. 

The plan of this article is as follows. In Section \ref{MellinVariables} the Cauchy problem satisfied by $\mathcal M_u(t, s)$, the Mellin transform of  suitable solutions $u$  of \eqref{eq:croisfrag2},\eqref{eq:croisfrag2data},\eqref{S1EK0}, is obtained. In Section \ref{Thm1} we prove Theorem~\ref{S1MainTh1} and Corollary \ref{S1Cor}.  In Section \ref{Extension} we study the extension of the local solution, and its uniqueness. The sign of the extension is studied in Section \ref{sign}, where Theorem \ref{S1MainTh4Bis} is proved.  Section \ref{SGammaNegatif} contains the case $\gamma <0$ and the proof of Theorem \ref{S1Th421}. Several  technical results are  gathered in the Appendix. The content of Sections \ref{MellinVariables} and \ref{Thm1}  where anounced and shortly presented in \cite{E}.

\section{The problem in Mellin variables} 
\label{MellinVariables}
We deduce in this Section the equation satisfied by the Mellin transform of a solution $u$ of \eqref{eq:croisfrag2} that would satisfy  suitable conditions. To this end we  suppose that $u(t, x)$ is a solution of equation \eqref{eq:croisfrag2} such that its  Mellin transform $\mathcal M_u$ is well defined for $s$ and $s+\gamma $, where $s$ belongs  to some domain $D$ of the complex plane $\C$.  Applying the  Mellin transform to both sides of equation \eqref{eq:croisfrag2} we arrive at:
\begin{eqnarray*}
\f{\p}{\p t} \mathcal M_u(t, s)+\int _0^\infty \f{\p}{\p x} \big(x^{1+\gamma } u (t,x)\big)x^{s-1}dx+\mathcal M_u(t, s+\gamma )=\int _0^\infty x^{s-1}\int_x^\infty \f{1}{y}k_0\left(\f{x}{y}\right) y^\gamma u(t,y)dydx\\
=\int _0^\infty dy y^{\gamma-1} u(t,y) \int_0^y  dx x^{s-1} k_0\left(\f{x}{y}\right)=\mathcal M_u(t, s+\gamma ) K(s).
\end{eqnarray*}
If $\lim _{ x\to 0 }x^{\gamma +s} u(t, x)=\lim _{ x\to \infty }x^{\gamma +s} u(t, x)=0$
we deduce that
\begin{eqnarray*}
\int _0^\infty \f{\p}{\p x} \big(x^{1+\gamma } u (t,x)\big)x^{s-1}dx=-(s-1)\int _0^\infty\big(x^{1+\gamma } u (t,x)\big)x^{s-2}dx=-(s-1)\mathcal M_u(t, s+\gamma )
\end{eqnarray*}
and finally,
\bqn
\label{S2eqmellin}
\f{\p}{\p t} \mathcal M_u(t, s)=( K(s)+s-2)\mathcal M_u(t, s+\gamma ).
\eqn
With our choice of  the measure $k_0$ (cf. \eqref{S1EK0} and \eqref{S1EKS0}), we are then led to consider the  problem 
\bqn
\frac {\partial W} {\partial t}(t, s)&=&\Phi (s)W(t, s+\gamma ),\,\,\,\forall s\in \C;\,\,\Re e(s)=s_*\label{S2eqconst}\\
W(0, s)&=&1,\,\,\forall s\in \mathscr{S}(\rho , \rho +\gamma ).\label{S2initial}
\eqn
for some $s_*\in \R$ and $\rho \in \R$, where
\bqn
\Phi (s)=\left( \frac {\theta} {s}+s-2\right),\,\,\forall s\in \C\setminus\{0\}. \label{S2eqconst2}
\eqn

Equations  like \eqref{S2eqconst}  have deserved some attention in the  literature, for a variety of functions $\Phi $ (cf.  \cite{BZ}, and references therein, \cite{EV1}) and have also been considered in \cite{BW}.  They may be Laplace transformed into a Carleman  type problem and solved using the classical Wiener-Hopf method (cf. \cite{C}, \cite{K}, \cite{PW}). See also Section \ref{Extension} for the same equation \eqref{S2eqconst} but a different  initial data.
\section{$\gamma >0$. Proof of Theorem \ref{S1MainTh1}}
\label{Thm1}
\subsection{An explicit solution of  \eqref{S2eqconst}-\eqref{S2initial}.}
\label{Explicit}
The  problem \eqref{S2eqconst}-\eqref{S2initial}  has a particularly simple and explicit solution:
\bqn
\label{S2solmellin}
\Omega (t, s)=F\left(\frac {s-\sigma _1} {\gamma }, \frac {s-\sigma _2} {\gamma }, \frac {s} {\gamma }, \gamma t\right)
\equiv  (1-\gamma t)^{\frac {2-s} {\gamma }}F\left(\frac {\sigma _1} {\gamma }, \frac {\sigma _2} {\gamma  }, \frac {s} {\gamma }, \gamma t \right).
\eqn
as it immediately follows from the identities  15.2.1  and  15.3.3 in \cite{AS}. We deduce from the properties of hypergeometric functions that,  for all $\rho >0$, $R>\rho $ and $T\in (0, \gamma ^{-1})$:
\bqn
&&\forall t\in (0, \gamma ^{-1}): \Omega (t, \cdot)\,\,\,\hbox{is analytic in the domain}\,\,\, s\in \C\setminus \{-m\gamma, m\in \N\} \label{S3EPO1}\\
&&\Omega  \in C\left([0, \gamma ^{-1})\times \mathscr{S}(0, \infty)\right)\label{S3EPO2}\\
&&\sup\left\{|\Omega (t, s)|,\,\,t\in [0, T),\,\,s\in \overline{\mathscr{S}(\rho, R)} \right\}<\infty\label{S3EPO3}\\
&&\lim _{ t\to 0 }\Omega (t, s)=1,\,\,\,\forall s\in \C\setminus \{-m\gamma , \,\,m\in \N\}\label{S3EPO4}\\
&&\lim _{ t\to \gamma ^{-1} }\Omega (t, s)=\Omega (\gamma ^{-1}, s)=
\frac {\Gamma \left( \frac {s} {\gamma  }\right)\Gamma \left( \frac {2-s} {\gamma  }\right)} {\Gamma \left( \frac {\sigma _1} {\gamma  }\right)\Gamma \left( \frac {\sigma _2} {\gamma  }\right)},\,\,\,\forall s\in \mathscr{S}(-\infty, 2)\setminus\{-n\gamma, n\in \N\}.\label{S3EPO5}
\eqn
Our purpose is now to take the inverse Mellin transform of the function $\Omega (t, s)$. We first notice:

\begin{proposition}
\label{S3P10}
For any $\sigma _0>0$, $\sigma _1$, $\sigma _2$   and $x >0$ and $T\in \left(0, \frac {1} {\gamma }\right)$ there exists a positive constant
$C=C(T, \sigma _0, \sigma _1, \sigma _2)$ such that:
$$
\int  _{ \Re e s =\sigma _0 }\left|\Omega (t, s)-(1-\gamma t)^{\frac {2-s} {\gamma }}\left(1+\frac {2 t} {s} \right)x^{-s}\right|ds\le 
C\left( x(1-\gamma t)^{\frac {1} {\gamma }}\right)^{-\sigma _0},\,\,\,\forall t\in (0, T).
$$
\end{proposition}
\proof  Using the series representation of  the Hypergeometric function in \eqref{S2solmellin}  we deduce, for each  $t\in (0, \gamma ^{-1})$ fixed:
$$
F\left(\frac {\sigma _1} {\gamma }, \frac {\sigma _2} {\gamma }, \frac {s} {\gamma }, \gamma t\right)-
1-\frac {2 t } {s }=\mathcal O _{ t } \left(|s|^{-2} \right),\,\,\,|s |\to \infty,\,\,\Re e s=\sigma _0>0,
$$
(cf. formula 15.2.2 in \cite{O}) and  there exists a constant $C=C(T, \sigma _0, \sigma _1, \sigma _2)>0$ such that if $t\in (0, T)$ and $s=\sigma _0+iv$, $v\in \R$ and $-\sigma _0\not \in \N$:
\bqn
\left|F\left(\frac {\sigma _1} {\gamma }, \frac {\sigma _2} {\gamma }, \frac {s} {\gamma }, \gamma t\right)-1-\frac {2 t } {s }\right|\le C(1+|s|)^{-2}.
\eqn
Then, by  definition of $\Omega (t, s)$:
\bqns
\left|\Omega (t, s)-(1-\gamma t)^{\frac {2-s} {\gamma }} \left( 1+\frac {2t} {s}\right)\right|\le C(1-\gamma t)^{\frac {2-s} {\gamma }}  (1+|s|)^{-2}
\eqns
and
\bqns
\int  _{ \Re e s =\sigma _0 }\left|\Omega (t, s)-(1-\gamma t)^{\frac {2-s} {\gamma }}\left(1+\frac {2 t} {s} \right)x^{-s}\right|ds\le 
C(1-\gamma t)^{\frac {2} {\gamma }}  \int  _{ \Re e s =\sigma _0 }\hskip -0.3cm(1+|s|)^{-3}\left|x(1-\gamma t)^{\frac {1} {\gamma }}\right|^{-s}ds\\
\le C(1-\gamma t)^{\frac {2} {\gamma }}  \left(x(1-\gamma t)^{\frac {1} {\gamma }}\right)^{-\sigma _0}
\int  _{ \Re e s =\sigma _0 }\hskip -0.3cm(1+|s|)^{-2}ds.
\eqns
\qed
It follows from Proposition \ref{S3P10} that  $\Omega (t, s)$ has an inverse Mellin transform when $0<\gamma t<1$. Our next purpose is to obtain its explicit expression.

\subsection{The inverse Mellin transform of $\Omega (t, s)$.}
\label{S4}
We recall that, for suitable functions $V$, the  classical  inverse Mellin transform  is defined as
\bqn
\label{S5EMT}
\mathcal M _{ \sigma _0 }^{-1}(V)=\frac {1} {2i\pi }\int  _{ \Re e(s)=\sigma _0 }x^{-s}V(s)ds
\eqn
for some $\sigma_0>0$ fixed. We  first  show the following:

\begin{proposition}
\label{S4P154}
Suppose $\sigma _1\in \C$, $\sigma _2\in \C$,  $\gamma >0$ and $t>0$ such that $0<\gamma t<1$ and define the function
\bqn
\label{S4EV}
v(t, x)=F\left(1+\frac {\sigma _1} {\gamma }, 1+\frac {\sigma _2} {\gamma }, 2, \gamma t\left(1+(\gamma t-1)x^\gamma  \right) \right)H\left(1-(1-\gamma t)^{\frac {1} {\gamma }}x\right)
\eqn
for $x>0$, where $H$ is the Heaviside function. Then, for all $s>0$, $v(s)\in E' _{ 0, q }$ for all $q>0$,  the Mellin transform of $v$ is given by:
\bqn
\label{S2Mellinv}
\mathcal M_v(t, s)\equiv \int _0^\infty v(t, x)x^{s-1}dx= (1-\gamma t)^{-\frac {s} {\gamma }}\frac {F\left(\frac {\sigma _1} {\gamma }, \frac {\sigma _2} {\gamma }, \frac {s} {\gamma }, \gamma t\right)-1} {\sigma _1\sigma _2t}.
\eqn
\end{proposition}

\proof For all $t$ fixed, $v(t)$ is a function with compact support in $(0, (1-\gamma t )^{-1/\gamma }))$, integrable on $(0, \infty)$ and bounded as $x\to 0$. Therefore $v\in \mathscr{C}\left([0, \gamma ^{-1}); E'_{0, q} \right)$ for all $q>0$. The proof of \eqref{S2Mellinv} is a straightforward calculation using the expression of the hypergeometric function. Since $\gamma >0$ we have that $\gamma t>0$. Moreover, since $1-\gamma t>0$ and $x>0$, we have $(\gamma t-1)x^\gamma <0$ and then $\left(1+(\gamma t-1)x^\gamma  \right)<1$. Finally, due to the  Heaviside function we only have to consider values of $(t, x)$ where $1+(\gamma t-1)x^\gamma >0$. It then follows from the expression of the function $F\left(1+\frac {\sigma _1} {\gamma }, 1+\frac {\sigma _2} {\gamma }, 2, \gamma t\left(1+(\gamma t-1)x^\gamma  \right) \right)$ as an absolutely convergent series (cf. definition 15.1.1 in \cite{AS}):
\bqns
\mathcal M _{ v }(t, s)=\sum_{ n=0 } ^\infty \frac {\Gamma \left(1+\frac {\sigma _1} {\gamma }+n \right)\Gamma \left( 1+\frac {\sigma _2} {\gamma }+n\right)\Gamma \left( 2\right)(\gamma t)^n} {\Gamma \left( 1+\frac {\sigma _1} {\gamma }\right)\Gamma \left( 1+\frac {\sigma _2} {\gamma }\right)\Gamma \left(2+n \right)\Gamma \left(n+1 \right)} 
\int_0^{(1-\gamma t)^{-\frac {1} {\gamma }}}\hskip -0.8cm(1+(\gamma t-1)x^\gamma )^n x^{s-1}dx.
\eqns
We use now that, since $\gamma >0$, we have for all $s>0$:
\bqns
\int_0^{(1-\gamma t)^{-\frac {1} {\gamma }}}(1+(\gamma t-1)x^\gamma )^n x^{s-1}dx&=&
\sum_{ m=0 } ^n \begin{pmatrix}n \\ m \end{pmatrix}(-1)^m(1-\gamma t)^m \int_0^{(1-\gamma t)^{-\frac {1} {\gamma }}}x^{m\gamma +s-1}dx= \label{S4E120}\\
&=&(1-\gamma t)^{-\frac {s} {\gamma }}\sum_{ m=0 } ^n \begin{pmatrix}n \\ m \end{pmatrix}\frac {(-1)^m} {m \gamma +s}=(1-\gamma t)^{-\frac {s} {\gamma }}\frac {\Gamma (n+1)\Gamma \left( \frac {s} {\gamma }\right)} {\gamma \Gamma \left(1+\frac {s} {\gamma }+n \right)}.\nonumber
\eqns

Then,
\bqns
\mathcal M_v(t, s)&=&(1-\gamma t)^{-\frac {s} {\gamma }}\sum_{ n=0 } ^\infty \frac {\Gamma \left(1+\frac {\sigma _1} {\gamma }+n \right)\Gamma \left( 1+\frac {\sigma _2} {\gamma }+n\right)\Gamma \left( 2\right)(\gamma t)^n} {\Gamma \left( 1+\frac {\sigma _1} {\gamma }\right)\Gamma \left( 1+\frac {\sigma _2} {\gamma }\right)\Gamma \left(2+n \right)\Gamma \left(n+1 \right)}\frac {\Gamma (n+1)\Gamma \left( \frac {s} {\gamma }\right)} {\gamma \Gamma \left(1+\frac {s} {\gamma }+n \right)}\\
&=& (1-\gamma t)^{-\frac {s} {\gamma }}\frac {F\left(\frac {\sigma _1} {\gamma }, \frac {\sigma _2} {\gamma }, \frac {s} {\gamma }, \gamma t\right)-1} {\sigma _1\sigma _2t}.
\eqns 
and this proves  \eqref{S2Mellinv}.
\qed
\vskip 0.3cm 
The next Corollary follows  from Proposition \ref{S4P154} and Theorem  11.10.1 in \cite{ML} on the uniqueness of the inverse Mellin transform:
\begin{corollary}
\label{coro1}
For all $\sigma _1\in \C$, $\sigma _2\in \C$,  suppose that $\gamma >0$, $0<\gamma t<1$ and let $u $ be the measure:
\bqns
u (t, x)=(1-\gamma t)^{\frac {1} {\gamma }}\delta \left(x-(1-\gamma t)^{-\frac {1} {\gamma }} \right)+ \sigma_1\sigma _2 t(1-\gamma t)^{\frac {2} {\gamma }}v(t, x).
\eqns
Then, for all $t\in \left(0, \gamma ^{-1}\right)$:
\bqns
\mathcal M_u(t)=\Omega (t),\,\,\,\hbox{and}\,\,\,
u(t)=\mathcal M^{-1}(\Omega (t)).
\eqns
\end{corollary}

We prove now the existence part in Theorem \ref{S1MainTh1}.
\vskip 0.15cm
\noindent
\begin{proposition}
\label{S3P10}
The measure $u$ defined in \eqref{S1Esol0}, \eqref{S1Esol} is a weak non negative solution of 
 \eqref{eq:croisfrag2},\eqref{S1EK0} on $(0, \gamma ^{-1})$ such that
\bqns
&&(i) \quad u \in \mathscr{C}\left(\left[0, \gamma ^{-1}\right); \mathscr{M} _{ \rho-1  } \right),\,\,\forall \rho >0. \\
&&(ii) \quad  \forall T\in \left(0, \gamma ^{-1}\right),\,\,\exists C_T>0;\,\, \int _0^\infty u(t, x)x^{s-1} dx\le C_T\,\,\,\,\forall t\in [0, T], \,\,\forall s >0
\eqns
and  satisfies  \eqref{S1E42}.
\end{proposition}
\textbf{Proof of Proposition \ref{S3P10}.} The assertions (i) and (ii) follows from the explicit expression of $u$. It is easy to  check that \eqref{S1E42} holds true.
Let us prove that $u$ is a weak solution of  \eqref{eq:croisfrag2},\eqref{S1EK0} on $(0, \gamma ^{-1})$.

We already  know that $\Omega (t, s)$ solves \eqref{S2eqconst} for all $t\in (0, \gamma ^{-1})$ and all $s\in \mathscr{S}(0, \infty)$. Since $\Omega (t)$ and $\Phi \Omega(t) $ are analytic and bounded in 
$\mathscr{S}(0, q)$ for any $q>0$, we deduce from Theorem 11.10.1 in \cite{ML} that $u\in \mathscr{C}((0, \gamma ^{-1}), E' _{ 0, q })$ for all $q>0$.
Applying the inverse Mellin tranform \eqref{S5EMT} to both sides of the equation  \eqref{S2eqconst} we deduce the following identity:
\bqn
\label{S4E99}
\frac {\partial u} {\partial t} (t, x)=\mathcal M _{ s_0 }^{-1}\left(\left( \frac {\theta} {s}+(s-1)-1\right)\tau  _{ \gamma  }\mathcal M_u \right)(t, x)
\eqn
where all the terms are in  $\mathscr{C}((0, \gamma ^{-1}), E' _{ 0, q })$ and we have denoted $(\tau  _{ \gamma  }\mathcal M_u)(t, s)=\mathcal M_u(t, s+\gamma )$. 
We consider now each of the terms in the right and side separately.  
Since $\sigma _0>0$, $\gamma >0$, using that $\mathcal M_u (t, s )=\Omega (t, s)$ for all $\Re e(s)>0$ we have:
\bqn
&&\mathcal M _{ \sigma _0 }^{-1}\left(\tau  _{ \gamma  }\mathcal M _{ u }\right)=x^\gamma u(t, x).\label{S4E100}\\
&&\mathcal M _{ s_0 }^{-1}\left((s-1)\tau  _{ \gamma  }\mathcal M_u\right)(t, x)=-\frac {\partial } {\partial x}\left(x^{\gamma +1}u(t, x)\right)\label{S4E101}
\eqn
In  the last term in  the right hand side of (\ref{S4E99}) we write as above:
$$
\frac {1} {2i\pi }\int  _{ \Re e s=\sigma _0 }\frac {\theta} {s}\mathcal M_u (t, s+\gamma )x^{-s}ds=
\int _0^\infty u(t, y)\left(\frac {1} {2i\pi }\int  _{ \Re e s=\sigma _0 }\frac {\theta} {s}y^{s+\gamma -1}x^{-s}ds\right)dy.
$$
Using that for $\sigma _0>0$:

\bqn
\label{S3E1650}
\frac {1} {2i\pi }\int  _{ \Re e s=\sigma _0 }\frac {1} {s}y^{s+\gamma -1}x^{-s}ds=
\begin{dcases}
0, & \hbox{if}\,\,y<x \\
 y^{\gamma -1}, &  \hbox{if}\,\, y>x
\end{dcases}
\eqn
we deduce
\bqn
\label{S4E102}
\frac {1} {2i\pi }\int  _{ \Re e s=\sigma _0 }\frac {\theta} {s}\mathcal M_u (t, s+\gamma )x^{-s}ds=
\theta \int _x^\infty u(t, y)y^{\gamma -1}dy.
\eqn
The left and right hand sides of equation \eqref{eq:croisfrag2},\eqref{S1EK0}  are then equal in $ \mathscr{C}((0, \gamma ^{-1}), E' _{ 0, q })$ and  in particular, $u\in\mathscr{C}((0, \gamma ^{-1}), \mathscr{D}'  (0, \infty))$ and  is a weak solution of  \eqref{eq:croisfrag2},\eqref{S1EK0}.

In order to prove the  non negativity of the measure $u$ we  use its definition  \eqref{S1Esol0}-\eqref{S1Esol} and the expression of the hypergeometric function in the right hand side of \eqref{S1Esol} as an absolutely convergent series:
$$
F\left(1+\frac {\sigma _1} {\gamma }, 1+\frac {\sigma _2} {\gamma }, 2, z \right)=
\sum_{ n=0 }^\infty \frac {\Gamma \left(1+\frac {\sigma _1} {\gamma }+n\right)\Gamma \left(1+\frac {\sigma _2} {\gamma }+n\right)z^n} 
{\Gamma \left(1+\frac {\sigma _1} {\gamma }\right)\Gamma \left(1+\frac {\sigma _2} {\gamma }\right)\Gamma (2+n)\Gamma (n+1)} 
$$
where we have denoted $\gamma  t\left(1+\left(\gamma  t-1\right)x^\gamma  \right)=z$.
When $\theta\in (0, 1)$ all the terms of the series are obviously non negative since $\sigma _2>0$ and $\sigma _1>0$. When $\theta>1$, we use that, since  $\sigma_1=\overline{\sigma _2}$, 
$\Gamma \left(1+\frac {\sigma _1} {\gamma }+n\right)=\overline{\Gamma \left(1+\frac {\sigma _2} {\gamma }+n\right)}$ for all  $n\in \N$, and again all the term of the series are non negative.
\qed
\vskip 0.4cm 
\begin{remark} 
The particular form of the measure $u(t)$ and  a simple calculation with  distributions in $(0, \infty)$ shows that the measure $u^S$ solves:
\bqns
\frac {\partial u^S(t)} {\partial t}+\frac {\partial } {\partial x}\left(x^{\gamma +1}u^S(t)\right)+x^\gamma u^S(t)=0,\,\,\,\hbox{in}\,\,\,\mathscr{D}_1'\left(\left(0, \gamma ^{-1}\right)\times (0, \infty)\right)
\eqns
and the function $u^R$ satisfies, for  all $t\in (0, \gamma ^{-1})$ and $x\in \left(0, (1-\gamma t)^{-1/\gamma }\right)$:
\bqns
\frac {\partial u^R(t, x)}{\partial t}+\frac {\partial \left(x^{\gamma +1}u^R(t, x)\right)} {\partial x}+x^\gamma u^R(t, x)=\theta\int _x^{(1-\gamma  t)^{-\frac {1} {\gamma }}}
\hskip -0.4cm u^R(t, y)y^{\gamma -1}dy+\theta \left(1-\gamma  t\right)^{-1}.
\eqns
\end{remark}
\vskip 0.15cm
\begin{remark}
\label{S3E564}
By the particular form of $u$ we deduce that, for all $\varphi \in C_c^1([0, \gamma ^{-1}]\times (0, \infty)):$
\bqn
\label{S3E565}
&&\int _0^\infty\int _0^\infty\left(\frac {\partial \varphi } {\partial t}+x^{\gamma +1}\frac {\partial \varphi } {\partial  x}+x^\gamma \varphi  (t, x) \right)u(t, x)dxdt
+\varphi (0, 1)=\nonumber\\
&&\hskip 0.3cm =\int _0^\infty u\left(\gamma ^{-1}, x \right)\varphi \left(\gamma ^{-1}, x \right)dx -\theta \int _0^\infty\int _0^\infty u(t, y)y^{\gamma -1} \int _0^y \varphi (t, x)dxdydt \label{S1Esol21}
\eqn
\end{remark}
\noindent
\textbf{Proof of Theorem \ref{S1MainTh1}.} By the Proposition \ref{S3P10} only the uniqueness of non negative  weak solutions satisfying \eqref{S1E40987}-\eqref{S1E42} for some $\rho >0$ remains to be proved. Suppose on the contrary  that $u$ and $v$ are two such solutions and $u(t)\not = v(t)$. Since $u_1$ and $u_2$ are nonnegative and satisfy \eqref{S1E40987} it follows that their Mellin transforms  $\mathcal M _{ u}(t)$ and  $\mathcal M _{v }(t) $  are well defined and analytic on $\mathscr{S}(\rho-\delta , \rho +\gamma +\delta  )$ for all $t\in (0, \gamma ^{-1})$ and satisfy \eqref{S4ES10}. By \eqref{S1E40987}, $\mathcal M _{ u}(t)$ and  $\mathcal M _{v }(t) $ also satisfy \eqref{S4ES678}. We check now that $u$ and $v$ also  satisfy \eqref{S4ES6780}.  Since the proof is of course  the same for both, we only consider $u$. By \eqref{S1E40987} and the continuity of the Mellin transform on $E' _{ \rho , \rho +\gamma  }$, it follows that  $\mathcal M _{ u }(0, s)=1$ for all $s\in\mathscr{S}(\rho, \rho +\gamma   )$ and $u$ satisfies \eqref{S4ES6780}.
Therefore, $M _{ u }$ and $M _{ v }$ satisfy the hypothesis of Theorem \ref{S4TU} and are then equal. This contradicts our hypothesis that  $u_1(t)\not = u_2(t)$, and proves the uniqueness. Assertion \eqref{S1Esol1000} has been shown in Proposition \ref{S3P10}.
\qed
\vskip 0.15cm
\noindent
\begin{remark}
\label{S3Run}
The existence of a non negative and non identically zero solution for the equation \eqref{eq:croisfrag2} with zero initial data  has been proved in \cite{BW} for quite general dislocation measures $k_0$ under some conditions. One of these conditions, denoted $(M_-)$, requires to have  $\sigma _1-1>0$.  That is not possible in our case by our choice of $k_0$ and \eqref{S1EKS024}.

\end{remark}
\vskip 0.15cm
\noindent
\textbf{Proof of Corollary \ref{S1Cor}}.
The behavior of $u(t, x)$ as $t\to \gamma ^{-1}$ is given by that of $F\left(1+\frac {\sigma _1} {\gamma }, 1+\frac {\sigma _2} {\gamma }, 2, z \right)$ as $z\to 1^-$ and depends on the values of $\frac {\sigma  _1} {\gamma }$ and $\frac {\sigma _2} {\gamma }$. 
\bqns
&&\lim_{ \gamma  t\to 1^- }u^R(t, x)=\lim _{ \gamma  t \to 1^- }\sigma _1\sigma _2t \left(1-\gamma  t\right)^{\frac {2} {\gamma }} F\left(1+\frac {\sigma _1} {\gamma }, 1+\frac {\sigma _2} {\gamma }, 2, \gamma  t\left(1+\left(\gamma  t-1\right)x^\gamma  \right) \right)
\\
&&=\lim_{ \gamma  t\to 1^- } \frac {\sigma _1\sigma _2t \left(1-\gamma  t\right)^{\frac {2} {\gamma }} } {\left(1-\gamma  t\left(1+\left(\gamma  t-1\right)x^\gamma  \right)\right)^{\frac {2} {\gamma }}}\lim _{ \gamma  t\to 1^- } \frac { F\left(1+\frac {\sigma _1} {\gamma }, 1+\frac {\sigma _2} {\gamma }, 2, \gamma  t\left(1+\left(\gamma  t-1\right)x^\gamma  \right) \right)} {\left(1-\gamma  t\left(1+\left(\gamma  t-1\right)x^\gamma  \right)\right)^{-\frac {2} {\gamma }}}
\eqns
Since $\frac {\sigma _1+\sigma _2} {\gamma }=\frac {2} {\gamma }>0$,  we have by 15.4.23 in \cite{O}:
\bqns
 \lim _{ \gamma  t\to 1^- } \frac { F\left(1+\frac {\sigma _1} {\gamma }, 1+\frac {\sigma _2} {\gamma }, 2, \gamma  t\left(1+\left(\gamma  t-1\right)x^\gamma  \right) \right)} {\left(1-\gamma  t\left(1+\left(\gamma  t-1\right)x^\gamma  \right)\right)^{-\frac {2} {\gamma }}}= \frac {\Gamma \left(\frac {2} {\gamma } \right)}
{ \Gamma \left(1+\frac {\sigma _1} {\gamma }  \right) \Gamma \left(1+\frac {\sigma _2} {\gamma }  \right)},\,\,\forall x>0.
\eqns
Since $1-\gamma  t\left(1+\left(\gamma  t-1\right)x^\gamma  \right)=(1-\gamma  t)(1+\gamma  t x^\gamma )$
it follows:
\bqns
\lim_{ \gamma  t\to 1^- } \frac {\sigma _1\sigma _2t \left(1-\gamma  t\right)^{\frac {2} {\gamma }} } {\left(1-\gamma  t\left(1+\left(\gamma  t-1\right)x^\gamma  \right)\right)^{\frac {2} {\gamma }}}
=\lim_{ \gamma  t\to 1^- } \frac {\sigma _1\sigma _2t  } 
{\left(1+\gamma  t x^\gamma \right)^{\frac {2} {\gamma }}}=
\frac {\sigma _1\sigma _2 } 
{\gamma \left(1+x^\gamma \right)^{\frac {2} {\gamma }}}
\eqns
and finally,
$$
\lim _{ \gamma  t\to 1^- }u^R(t, x)= \frac {\Gamma \left(\frac {2} {\gamma } \right)\sigma _1\sigma _2}
{\gamma  \Gamma \left(1+\frac {\sigma _1} {\gamma }  \right) \Gamma \left(1+\frac {\sigma _2} {\gamma }  \right)}(1+x^\gamma )^{-\frac {2} {\gamma }}=\frac {\gamma \Gamma \left(\frac {2} {\gamma } \right)}
{  \Gamma \left(\frac {\sigma _1} {\gamma }  \right) \Gamma \left(\frac {\sigma _2} {\gamma }  \right)}(1+x^\gamma )^{-\frac {2} {\gamma }}.
$$
The proof of properties (\ref{S1Ebw1})--(\ref{S1Ebw3}) follow from the explicit expression  (\ref{S2solmellin}) of $\Omega (t, s)$  and  formulas 15.4(ii) in \cite{O}.
\qed

\section{$\gamma >0$. Extension of the local solution.}
\label{Extension}
The two main results of this Section  are the following. The first, where we use the notation
\bqn
\label{S1Enu}
\nu=\min (2, \Re e(\sigma _2)+\gamma).
\eqn
extends the local solution $u$ obtained in Theorem \ref{S1MainTh1}:
\begin{theorem}
\label{S1MainTh3}  
For all $\theta>0$ and $\gamma >0$ there exists a  global weak solution  $w \in  \mathscr{C}\left([0, \infty); E' _{ 0, \nu }\right)$  of  \eqref{eq:croisfrag2},\eqref{S1EK0} on $t\in (0, \infty)$,  such that $w(t)=u(t)$ for all $t\in (0, \gamma ^{-1})$ and satisfying
\bqn
\label{S1E198765}
\left\{
\begin{split}
&w \in \mathscr{C}^\infty((\gamma ^{-1},\infty)\times (0, \infty))\cap  \mathscr{C}([\gamma ^{-1},\infty)\times (0, \infty))\\
&w\,\,\, \hbox{satisfies}\,\,  \eqref{eq:croisfrag2},\eqref{S1EK0},\,\, \hbox{pointwise},\,\,\forall t>\gamma ^{-1}, \forall x>0\\
&\forall t>\gamma ^{-1},\,\,\mathcal M_w(t, \cdot)  \,\,\hbox{is analytic in} \,\,\mathscr{S}(0, \Re e(\sigma _2)+\gamma ),\\
&\mathcal M_w (\gamma ^{-1},s)=
\Omega (\gamma ^{-1}, s)\,\,\, \forall s\in \C,\,\,\,\Re e(s)<2
\end{split}
\right.
\eqn
\end{theorem}
\vskip 0.15cm
\noindent
An  expression of $w $ for $t>\gamma ^{-1}$ is  obtained in Remark \ref{S4Eexpl}. 

Then, the following uniqueness result is proved, for  the case $\gamma \in (0, 2)$:
\begin{theorem}
\label{S1MainTh2BIS} If $\gamma \in (0, 2)$, the measure $w$ defined in Theorem \ref{S1MainTh3}
is the unique global weak solution of  \eqref{eq:croisfrag2},\eqref{S1EK0} such that for all $T>0$,  $w \in  \mathscr{C}\left([0, T); E' _{ \rho -\delta , \rho +\gamma +\delta }\right)$  for some $\rho\in (0, 2-\gamma )$, $\delta >0$ and 
\bqn
\label{S1E2400Bis}
\left\{
\begin{split}
&\sup\left\{|\mathcal M_w(s, t)|;\,\, s\in \mathscr{S}(\rho-\delta  , \rho +\gamma +\delta ),\,\,t\in \left[0, T\right] \right\}<\infty \label{S1E2402Bis}\\
&w(0)=\delta_1.
\end{split}
\right.
\eqn
 \end{theorem}
The uniqueness of the solution obtained in Theorem \ref{S1MainTh3}, valid for all $\gamma >0$, is  proved in Theorem \ref{S1MainTh2}.

In order to extend the solution $u$ of \eqref{eq:croisfrag2},\eqref{eq:croisfrag2data} beyond $t=\gamma ^{-1}$ we first obtain a solution of  \eqref{S2eqconst} for $\gamma t>1$.
\subsection{Another explicit solution of \eqref{S2eqconst}.}
\label{S6}
The new solution of the  equation \eqref{S2eqconst} is obtained  as follows:
\begin{proposition} 
\label{S5PropU}The function $U$ defined, for all $t>\gamma ^{-1}$ and all $s\in \C$, as 
\bqn
\label{S5EsolU}
U(t, s)=\frac{(\gamma  t)^{\frac {\sigma _1-s} {\gamma }}\Gamma \left( \frac {s} {\gamma }\right)\Gamma \left(1- \frac {s-\sigma _2} {\gamma }\right)}{\Gamma \left( \frac {\sigma _1} {\gamma }\right)\Gamma \left(1- \frac {\sigma _1-\sigma _2} {\gamma }\right)}
F\left( 1-\frac {\sigma _1} {\gamma }, \frac {s-\sigma _1} {\gamma }; 1+\frac {\sigma _2-\sigma _1} {\gamma }; \frac {1} {\gamma t}\right)
\eqn
is such that:
\vskip 0.15cm
\noindent
(i) $U$ is meromorphic on $(\gamma ^{-1}, \infty)\times S$ and the set $S$ of its poles is:
\bqns
S=\left\{s=-m\gamma,\,\,\,m\in \N\right\}\cup \left\{s=\Re e(\sigma _2)+(m+1)\gamma , \,\,\,m\in \N\right\}.
\eqns
\vskip 0.15cm
\noindent
(ii) $U$ satisfies the equation \eqref{S2eqconst} for all $t>\gamma ^{-1}$ and all $s\in \C \setminus S$;
\vskip 0.15cm
\noindent
(iii) For all closed  subinterval $I\subset (0, \Re e(\sigma _2)+ \gamma ))$, there exists a positive constant $C=C(t, I)$ such that, for all $t\ge\gamma ^{-1}$ and all $s\in \C$, $\Re e(s)\in I$
\bqns
|U(t, s)|+|U_t(t, s)|\le Ce^{-\frac {\pi |\Im m(s)|} {\gamma }}t^{\frac {\sigma_1 -s_0} {\gamma }}\left((1+|s|)^{-1+\frac {\sigma _1} {\gamma }}+(1+|s|)^{-\frac {\sigma _2} {\gamma }} \right)
\eqns
\vskip 0.15cm
\noindent
(iv) for $\Re e(s) < 2:$\hskip 0.5cm  $\displaystyle{\lim _{ t\to \left(\gamma ^{-1}\right)^+ }U (t, s)=\lim _{ t\to \left(\gamma ^{-1}\right)^- }\Omega (t, s)}
=\frac{\Gamma \left( \frac {s} {\gamma }\right)}{\Gamma \left( \frac {\sigma _1} {\gamma }\right)}
\frac {\Gamma \left(\frac {\sigma _2+\sigma _1-s} {\gamma } \right)} 
{\Gamma \left(\frac {\sigma _2} {\gamma } \right)}$.

\end{proposition}
\textbf{Proof.}
It is straightforward to check that the function 
\bqn
\label{S5EV}
V(s)=\gamma ^{\frac {s} {\gamma }}\frac {\Gamma \left(\frac {s-\sigma _1} {\gamma } \right)} {\Gamma \left(\frac {s} {\gamma } \right)\Gamma \left(1-\frac {s-\sigma _2} {\gamma } \right)},\,\,\,\forall s\in \C\setminus \{s\in \C; s=\sigma _1-m\gamma ,\, m\in \N\}
\eqn
satisfies:
\bqn
\label{S5EqV}
V(s+\gamma )=-\frac {(s-\sigma _1)(s-\sigma _2)} {s}V(s),\,\,\,\forall s\in \C\setminus \{s\in \C; s=\sigma _1-m\gamma ,\, m\in \N\}
\eqn

We define now the function of $t$ and $s$:
\beq
\label{S5E146}
\mathscr{U}(t, s)=\frac {1} { \gamma V(s)}\int _{\mathscr{C}}\frac {(-t)^{\frac {\sigma -s} {\gamma }}V(\sigma )d\sigma } {\Gamma \left( 1+\frac {\sigma -s} {\gamma }\right)\left(e^{-\frac {2i\pi } {\gamma }(\sigma -s)}-1\right)}
\eeq
where the path of integration  $\mathscr{C}$ is  as follows. For $\sigma _0>\sigma _1$  fixed, we define:
\bqn
&&\mathscr{C}=\mathscr{C}_1\cup \mathscr{C}_2\cup \mathscr{C}_3\label{S5C}\\
&&\mathscr{C}_1=\left\{s\in \C; s=\sigma _0+iv, |v|\le 1\right\}\label{S5C1}\\
&&\mathscr{C}_2=\left\{s\in \C; s=(u+iv), \,\, v=-u+(1+\sigma _0),\,\,\, v>1\right\}\label{S5C2}\\
&&\mathscr{C}_3=\left\{s\in \C; s=(u+iv), \,\, v=u-(1+\sigma _0),\,\,\, v<-1\right\}.\label{S5C3}
\eqn
Notice that when $\sigma\in \mathscr{C}$ is such that $|\sigma |\to \infty$ we have that $\Re e(\sigma )\to -\infty$.  The function below the integral sign of \eqref{S5E146} may be written as follows:
\begin{eqnarray}
w(t, s, \sigma )=\frac {(-t)^{\frac {\sigma -s} {\gamma }}V(\sigma )} {\Gamma \left( 1+\frac {\sigma -s} {\gamma }\right)
\left(e^{-\frac {2i\pi } {\gamma }(\sigma -s)}-1\right)}=
\frac {(-t)^{-\frac { s} {\gamma }}(-\gamma t)^{\frac {\sigma} {\gamma }}\Gamma \left(\frac {\sigma -\sigma _1} {\gamma } \right)}
 {\Gamma \left( 1+\frac {\sigma -s} {\gamma }\right)
\left(e^{-\frac {2i\pi } {\gamma }(\sigma -s)}-1\right)\Gamma \left(\frac {\sigma } {\gamma } \right)\Gamma \left(1-\frac {\sigma -\sigma _2} {\gamma } \right)}
\label{S5E147}
\end{eqnarray}
For $s\in \C$ fixed,  $\sigma \in  \mathscr{C}$, $|\sigma |\to \infty$  we have, using Stirling's formula:
\bqns
\left|\Gamma \left(\frac {\sigma -\sigma _1} {\gamma } \right)\right|&\approx &e^{\frac {\sigma -\sigma _1} {\gamma }\log\left(\frac {\sigma -\sigma _1} {\gamma } \right)}
\approx e^{\frac {\Re e(\sigma )} {\gamma }\log |\sigma /\gamma |}\\
\left|\Gamma \left(\frac {\sigma } {\gamma } \right)\right| &\approx &
e^{\frac {\sigma} {\gamma }\log\left(\frac {\sigma} {\gamma } \right)}
\approx e^{\frac {\Re e(\sigma )} {\gamma }\log |\sigma /\gamma |}\\
\left|\Gamma \left(1-\frac {\sigma -\sigma_2} {\gamma } \right)\right| &\approx &
e^{\left(1-\frac {\sigma -\sigma _2} {\gamma }\right)\log\left(1-\frac {\sigma -\sigma _2} {\gamma } \right)}
\approx e^{\left(1-\frac {\Re e(\sigma )} {\gamma }\right)\log |\sigma /\gamma |}\\
\left|\Gamma \left(1+\frac {\sigma -s} {\gamma } \right)\right| &\approx &
e^{\left(1+\frac {\sigma -s} {\gamma }\right)\log\left(1+\frac {\sigma -s} {\gamma } \right)}
\approx e^{\left(1+\frac {\Re e(\sigma )} {\gamma }\right)\log |\sigma /\gamma |}\\
\left|e^{-\frac {2i\pi } {\gamma }(\sigma -s)}-1\right|&=& \left|e^{-\frac {2i\pi } {\gamma }(\Re e(\sigma )+i\Im m(\sigma ) )}-1\right|\ge  
\left|e^{\frac {2\pi \Im m(\sigma )} {\gamma }}-1\right|\\
(-\gamma t)^{\frac {\sigma -s} {\gamma }}&=&e^{\frac {\sigma -s} {\gamma }\log (-\gamma t)}\approx e^{\frac {\Re e(\sigma )} {\gamma }\log(\gamma t)}.
\eqns
We deduce that, for each  $s\in \C$ fixed there exists a constant $C=C(s)>0$ such that for all $t>\gamma ^{-1}$ and all $\sigma \in  \mathscr{C}$:
\bqns
\left|w(t, s, \sigma )\right|\le C\frac {e^{\frac {\Re e(\sigma )} {\gamma }\log(\gamma t)}} {e^{\log |\sigma /\gamma |}\left|e^{\frac {2\pi v} {\gamma }}-1\right|}
\eqns
where the constant $C$ depends on $s$.  Since $\Re e(\sigma )\to -\infty$ as $|\sigma |\to \infty$ for $\sigma \in  \mathscr{C}$ the function $w(t, s, \cdot)$ is exponentially decaying in $\sigma $ for $t>\gamma ^{-1}$, and  the integral in the right hand side of \eqref{S5E146} is absolutely convergent. It follows that the function $\mathscr{U}(t, s)$ is well defined, continuous with respect to $t$ and analytic with respect to $s$ for all $t>\gamma ^{-1}$ and $s$ in
\bqn
\label{S6EDom}
D=\C\setminus\left( \{s\in \C;\,\,s=-m\gamma ,\,\,m\in \N\}\cup \{s\in \C;\,\,s=\sigma _2+(m+1)\gamma ,\,\,m\in \N\}\right).
\eqn
By the exponential decay of $w(t, s, \sigma )$ in $\sigma $ along   $\mathscr{C}$ and its regularity in time,  a simple calculation yields:
\bqn
\frac {\partial\mathscr{U}} {\partial t}(t, s)=\frac {-1} { \gamma V(s)}\int _{\mathscr{C}}\frac {(-t)^{\frac {\sigma -s} {\gamma }-1}\left(\frac {\sigma  -s} {\gamma }\right)V(\sigma )d\sigma } {\Gamma \left( 1+\frac {\sigma -s} {\gamma }\right)\left(e^{-\frac {2i\pi } {\gamma }(\sigma -s)}-1\right)} \label{S6EderiU1}\\
=\frac {-1} { \gamma V(s)}\int _{\mathscr{C}}\frac {(-t)^{\frac {\sigma -s} {\gamma }-1}V(\sigma )d\sigma } {\Gamma \left( \frac {\sigma -s} {\gamma }\right)\left(e^{-\frac {2i\pi } {\gamma }(\sigma -s)}-1\right)} \label{S6EderiU2}\\
=\frac {-1} {\gamma  V(s)}\int _{\mathscr{C}}\frac {(-t)^{\frac {\sigma -(s+\gamma )} {\gamma }}V(\sigma )d\sigma }
 {\Gamma \left(1+ \frac {\sigma -(s+\gamma )} {\gamma }\right)\left(e^{-\frac {2i\pi } {\gamma }(\sigma -(s+\gamma ))}-1\right)}. \label{S6EderiU3}
\eqn
By \eqref{S5EqV}:
$$
\frac {-1} {V(s)}=\frac {(s-\sigma _1)(s-\sigma _2)} {s}\frac {1} {V(s+\gamma )}.
$$
We deduce
$$
\frac {\partial \mathscr{U}} {\partial t}(t, s)=\frac {(s-\sigma _1)(s-\sigma _2)} {s}\mathscr{U}(t, s+\gamma ).
$$
and the function $\mathscr{U}(t, s)$  satisfies the equation \eqref{S2eqconst} for $t>\gamma ^{-1}$ and $s$ as in \eqref{S6EDom}. 

Our next step is to prove that $\mathscr{U}=U$,  using the residue's method. To this end,  we notice that for $s$ fixed as in \eqref{S6EDom}, the poles of the function $w(t, s, \sigma )$ to integrate are located at the following points:
\bqns
\left\{\sigma =\sigma _1-\gamma m,\,\,m\in \N\right\}\cup
\left\{\sigma =s+\gamma m, \,\,m\in \N\right\}
\eqns
For  values of $t$ such that $\gamma t>1$ we must use the residues at the  points $\sigma =\sigma _1-\gamma m$. We deform the integration contour, always in the region where $\Re s\to -\infty$. Since $\sigma _0>\Re e(\sigma _1)$:
\bqns
\mathscr{U}(t,s)=\frac {2i\pi } {V(s)\left(e^{-\frac {2i\pi } {\gamma }(\sigma _1-s)} -1\right)}\sum_{ m=0 }^\infty  \frac {(- t)^{\frac {\sigma _1-s} {\gamma }-m}\gamma ^{\frac {\sigma _1} {\gamma }-m} (-1)^m  } 
{\Gamma \left( \frac {\sigma _1} {\gamma }-m\right)\Gamma \left(1+ \frac {\sigma _1-s} {\gamma }-m\right)\Gamma \left(1- \frac {\sigma _1-\sigma _2} {\gamma }+m\right)
\Gamma (m+1)}\\
=\frac {2i\pi  (- t)^{\frac {\sigma _1-s} {\gamma }}\gamma ^{\frac {\sigma _1} {\gamma }}} {V(s)\left(e^{-\frac {2i\pi } {\gamma }(\sigma _1-s)} -1\right)}\sum_{ m=0 }^\infty  \frac {(\gamma t)^{-m}} 
{\Gamma \left( \frac {\sigma _1} {\gamma }-m\right)\Gamma \left(1+ \frac {\sigma _1-s} {\gamma }-m\right)\Gamma \left(1- \frac {\sigma _1-\sigma _2} {\gamma }+m\right)
\Gamma (m+1)}\\
=\frac {2i\pi  (- t)^{\frac {\sigma _1-s} {\gamma }}\gamma ^{\frac {\sigma _1} {\gamma }}F\left( 1-\frac {\sigma _1} {\gamma }, \frac {s-\sigma _1} {\gamma }; 1-\frac {\sigma _1-\sigma _2} {\gamma }; \frac {1} {\gamma t}\right)} 
{V(s)\left(e^{-\frac {2i\pi } {\gamma }(\sigma _1-s)} -1\right)\Gamma \left( \frac {\sigma _1} {\gamma }\right)\Gamma \left(1- \frac {s-\sigma _1} {\gamma }\right)\Gamma \left(1- \frac {\sigma _1-\sigma _2} {\gamma }\right)}.
\eqns
This may we writen: 
\bqns
&&\mathscr{U}(t, s)= \frac {2i\pi (-\gamma  t)^{\frac {\sigma _1-s} {\gamma }}\Gamma \left( \frac {s} {\gamma }\right)\Gamma \left(1- \frac {s-\sigma _2} {\gamma }\right)F\left( 1-\frac {\sigma _1} {\gamma }, \frac {s-\sigma _1} {\gamma }; 1-\frac {\sigma _1-\sigma _2} {\gamma }; \frac {1} {\gamma t}\right)} 
{\Gamma \left( \frac {\sigma _1} {\gamma }\right)\Gamma \left(1- \frac {\sigma _1-\sigma _2} {\gamma }\right)\Gamma \left(\frac {s-\sigma _1} {\gamma }\right)\Gamma \left(1- \frac {s-\sigma _1} {\gamma }\right)\left(e^{-\frac {2i\pi } {\gamma }(\sigma _1-s)} -1\right)}\\
\eqns
and using the identity $\Gamma (x)\Gamma (1-x)(e^{2i\pi x}-1)=2i\pi e^{i\pi x}$:
\bqn
\mathscr{U}(t, s)&=&\frac{(-\gamma  t)^{\frac {\sigma _1-s} {\gamma }}}{\Gamma \left( \frac {\sigma _1} {\gamma }\right)\Gamma \left(1- \frac {\sigma _1-\sigma _2} {\gamma }\right)}
\Gamma \left( \frac {s} {\gamma }\right)\Gamma \left(1- \frac {s-\sigma _2} {\gamma }\right)\times \nonumber\\
&& \times F\left( 1-\frac {\sigma _1} {\gamma }, \frac {s-\sigma _1} {\gamma }; 1-\frac {\sigma _1-\sigma _2} {\gamma }; \frac {1} {\gamma t}\right) 
e^{-\frac {i\pi } {\gamma }(s-\sigma _1)} \label{S5ETH1}
\eqn
from where,  using that $e^{-i\pi }=-1$ it follows that $\mathscr{U}(t, s)=U(t, s)$ for $\gamma t>1$ and $s\in \C$. Property (ii) immediately follows. In order to prove the property (iii)  we use again  Stirling's formula. For $s\in \C$ such that $\Re e(s)$ remains bounded and $|\Im m(s)|\to\infty$:
\bqns
\left|\Gamma \left(\frac {s } {\gamma } \right)\right|\approx e^{\frac {iv} {\gamma }i Arg (s /\gamma )}=e^{-\frac {v} {\gamma }Arg (s /\gamma )}\\
\left|\Gamma \left(1-\frac {s -\sigma_2} {\gamma } \right)\right|\approx e^{\frac {-iv} {\gamma }i Arg (-s /\gamma )}=e^{\frac {v} {\gamma }Arg (-s /\gamma )}
\eqns
Then,
\bqns
\left|\Gamma \left(\frac {s } {\gamma } \right)\Gamma \left(1-\frac {s -\sigma_2} {\gamma } \right)\right|
\approx e^{-\frac {v} {\gamma }Arg (s /\gamma )}e^{\frac {v} {\gamma }Arg (-s /\gamma )}.
\eqns
If $v\to \infty$, $Arg(s/\gamma )\to \pi /2$, $Arg(-s/\gamma )\to -\pi /2$ and then
\bqns
\left|\Gamma \left(\frac {s } {\gamma } \right)\Gamma \left(1-\frac {s -\sigma_2} {\gamma } \right)\right|
\approx e^{-\frac {v} {\gamma }\frac {\pi } {2}}e^{\frac {v} {\gamma }\frac {(-\pi )} {2}}\approx e^{-\pi v}.
\eqns
If $v\to -\infty$, $Arg(s/\gamma )\to -\pi /2$ and $Arg(-s/\gamma )\to \pi /2$ from where
\bqns
\left|\Gamma \left(\frac {s } {\gamma } \right)\Gamma \left(1-\frac {s -\sigma_2} {\gamma } \right)\right|
\approx e^{\frac {-v} {\gamma }\frac {(-\pi } {2})}e^{\frac {v} {\gamma }\frac {\pi } {2}}\approx e^{\pi v}.
\eqns
On the other hand, by  15.7.2 in \cite{AS}:
$$
\left|F\left( 1-\frac {\sigma _1} {\gamma }, \frac {s-\sigma _1} {\gamma }; 1+\frac {\sigma _2-\sigma _1} {\gamma }; \frac {1} {\gamma t}\right)\right|=
\mathcal O\left(|s|^{-1+\frac {\sigma _1} {\gamma }}+|s|^{-\frac {\sigma _2} {\gamma }} \right),\,\,\,|\Im m(s)|\to \infty.
$$
Using $|(\gamma t)^{\frac {\sigma_1 -s} {\gamma }}|=(\gamma t)^{\frac {\sigma_1 -s_0} {\gamma }}$  we deduce:
\bqns
|U(t, s)|& \sim &Ce^{-\frac {v} {\gamma }\frac {\pi } {2}}e^{\frac {v} {\gamma }\left(-\frac {\pi } {2}\right)}t^{\frac {\sigma_1 -s_0} {\gamma }}
\left|F\left( 1-\frac {\sigma _1} {\gamma }, \frac {s-\sigma _1} {\gamma }; 1+\frac {\sigma _2-\sigma _1} {\gamma }; \frac {1} {\gamma t}\right)\right|\\
\hskip 5cm &=&e^{-\frac {\pi v} {\gamma }}t^{\frac {\sigma_1 -s_0} {\gamma }}\mathcal O\left(|s|^{-1+\frac {\sigma _1} {\gamma }}+|s|^{-\frac {\sigma _2} {\gamma }} \right),\,\,\,v\to \infty
\\
|U(t, s)|& \sim &Ce^{-\frac {v} {\gamma }\left(-\frac {\pi } {2}\right)}e^{\frac {v} {\gamma }\frac {\pi } {2}}t^{\frac {\sigma_1 -s_0} {\gamma }}
\left|F\left( 1-\frac {\sigma _1} {\gamma }, \frac {s-\sigma _1} {\gamma }; 1+\frac {\sigma _2-\sigma _1} {\gamma }; \frac {1} {\gamma t}\right)\right|\\
\hskip 5cm &=&e^{\frac {\pi v} {\gamma }}t^{\frac {\sigma_1 -s_0} {\gamma }}\mathcal O\left(|s|^{-1+\frac {\sigma _1} {\gamma }}+|s|^{-\frac {\sigma _2} {\gamma }} \right),\,\,\,v\to -\infty
\eqns
and (iii) follows.

The property  (iv)  is  directly deduced from the definitions of $\Omega $ and $U$, \eqref{S2solmellin} and \eqref{S5EsolU} respectively, and the  application of 
identity 15.4.20 in \cite{O}  when $s<\Re e(\sigma _2+\sigma _1)\equiv 2$:
\bqns
\lim _{ \gamma t\to 1^+ }F\left( 1-\frac {\sigma _1} {\gamma }, \frac {s-\sigma _1} {\gamma }; 1+\frac {\sigma _2-\sigma _1} {\gamma }; \frac {1} {\gamma t}\right)
&=&
\frac {\Gamma \left(1+\frac {\sigma _2-\sigma _1} {\gamma } \right)\Gamma \left(\frac {\sigma _2+\sigma _1-s} {\gamma } \right)} 
{\Gamma \left(\frac {\sigma _2} {\gamma } \right)\Gamma \left(1+\frac {\sigma _2-s} {\gamma } \right)}\\
\lim _{ \gamma t\to 1^- }F\left(\frac {s-\sigma _1} {\gamma }, \frac {s-\sigma _2} {\gamma }, \frac {s} {\gamma }, \gamma t\right)
&=&
\frac{\Gamma \left( \frac {s} {\gamma }\right)}{\Gamma \left( \frac {\sigma _1} {\gamma }\right)}
\frac {\Gamma \left(\frac {\sigma _2+\sigma _1-s} {\gamma } \right)} 
{\Gamma \left(\frac {\sigma _2} {\gamma } \right)}.
\eqns
\qed
\vskip 0.2 cm
\begin{remark}
\label{S4R34}
The solution of  \eqref{S2eqconst} defined in \eqref{S2solmellin} may be obtained in the same way as the solution $U$ has been obtained in Proposition \ref{S5PropU}. It is enough to this end to start with the following function $\widetilde V(s)$:
\beq
\label{S5E2}
\widetilde V(s)=(-\gamma )^{\frac {s} {\gamma }}\frac {\Gamma \left(\frac {s-\sigma _1} {\gamma }\right)\Gamma \left(\frac {s-\sigma _2} {\gamma }\right)} {\Gamma \left(\frac {s} {\gamma }\right)}
\eeq
that also satisfies the equation  \eqref{S5EqV}. It may then be checked that:

\beq
\label{S5E3}
\Omega (t, s)=\frac {-1} { \gamma \widetilde V(s)}\int _{\Re e \sigma =\sigma _0}\frac {(-t)^{\frac {\sigma -s} {\gamma }}\widetilde V(\sigma )d\sigma } {\Gamma \left( 1+\frac {\sigma -s} {\gamma }\right)\left(e^{-\frac {2i\pi } {\gamma }(\sigma -s)}-1\right)}
\eeq
for any $\sigma _0>0$.

\end{remark}
\subsection{Inverse Mellin transform of the function $U$ in \eqref{S5EsolU}. }
By the Proposition \ref{S5PropU} it is  possible to apply to the function $U(t)$  the  inverse  Mellin transform defined as in \eqref{S5EMT} with  $s_0\in (0, \Re e(\sigma _2))$. We then  define:
\bqn
\label{S6E642}
\omega (t, x)=\mathcal M _{ s_0 }^{-1}(U(t))(x).
\eqn
\begin{proposition}
\label{S6Propo2}
For any $s_0\in (0, \Re e(\sigma _2))$ the function $\omega $ defined in  \eqref{S6E642} is such that:
\vskip 0.15cm
\noindent
(i) $\omega \in C^{\infty }\left( (\gamma ^{-1}, \infty)\times (0, \infty)\right)$
\vskip 0.15cm
\noindent
(ii) $\omega \in \mathscr{C}((\gamma ^{-1}, \infty); E' _{ 0, \Re e(\sigma _2)+\gamma })$ and $\mathcal M _{ \omega  }(t, s)=U(t, s)$,  for all $t>\gamma ^{-1}, s\in \mathscr{S}(0, \Re e(\sigma _2)+\gamma )$, 
\vskip 0.15cm
\noindent
(iii) $\omega $ satisfies the equation \eqref{eq:croisfrag2},\eqref{S1EK0} pointwise for all $t>\gamma ^{-1}$ and $x>0$.
\vskip 0.15cm
\noindent
(iv) For all $x>0$:
\bqns
\lim _{ t\to \left(\gamma ^{-1}\right)^+ }\omega (t, x)= \frac{\gamma\Gamma \left( \frac {2} {\gamma }\right)}{\Gamma \left( \frac {\sigma _1} {\gamma }\right)\Gamma \left( \frac {\sigma _2} {\gamma }\right)}
\left(1+x^\gamma  \right)^{-\frac {2} {\gamma }}
\eqns
(v) For all $t>\gamma^{-1}:$ 
\bqns
\omega (t, x)=\frac{\gamma \Gamma \left(1+ \frac {\sigma _2} {\gamma }\right)}
{\Gamma \left( \frac {\sigma _1} {\gamma }\right)\Gamma \left(1+ \frac {\sigma _2-\sigma _1} {\gamma }\right)}\frac {(\gamma t-1)^{\frac {\sigma _1} {\gamma }-1}} {(\gamma t)^{\frac {\sigma _2} {\gamma }}}x^{-\sigma _2-\gamma }+o\left(x^{-\sigma _2-\gamma}  \right),\,\,\,x\to \infty
\eqns
\vskip 0.15cm
\noindent
\end{proposition}
\textbf{Proof}
The assertion (i) It  follows from the regularity of the function $U$ with respect to $t$ and  assertion $(iii)$ in Proposition \ref{S5PropU}.

For the proof of  assertion $(ii)$  we first notice that, by properties (i) and (iii) of Proposition \ref{S5PropU}, the function $U(t, s)$ is analytic and exponentially decaying  the strip 
$s\in \mathscr{S}(0,  \Re e(\sigma _2)+\gamma )$ as $|\Im m(s)|\to \infty$. Then, by classical properties of the Mellin transform (cf. Theorem 11.10.1  in \cite{ML}) assertion (ii) follows.

In order to  prove assertion (iii) we first notice that,  from the analyticity and boundedness  properties of $\Phi U(t)$ on the strip $s\in \mathscr{S}(0, \Re e(\sigma _2)+\gamma )$ 
and by the same general properties of the Mellin transform,  we have $\Phi U\in \mathscr{C}((\gamma ^{-1}, \infty); E' _{ 0, \Re e(\sigma _2)+\gamma  })$.
We may then apply the inverse Mellin transform as defined in \eqref{S5EMU2} to both sides of the equation \eqref{S2eqconst}

\bqn
\label{S6E345}
\frac {\partial } {\partial t}\mathcal M _{ s_0 }^{-1}(U(t))(x)=\mathcal M _{ s_0 }^{-1}\left(\Phi  U   (t, \cdot+\gamma )\right)(x).
\eqn
By the the regularity of $U(t, s)$ with respect to $t$ and the decay properties of $U(t, s)$ and $U_t(t, s)$  along the integration curve $\Re e(s)=s_0$ we have
\bqn
\label{S2E789}
\frac {\partial } {\partial t}\mathcal M _{ s_0 }^{-1}(U(t))(x)=\frac {\partial \omega } {\partial t}(t, x). 
\eqn
Arguing now as in the proof of Proposition \ref{S3P10} we deduce that $\omega $ satisfies the equation \eqref{eq:croisfrag2},\eqref{S1EK0}  where all  the  terms belong to 
$\mathscr{C}((\gamma ^{-1}, \infty); E' _{ 0, \Re e(\sigma _2)+\gamma  })$. Therefore, the equation is satisfied in the weak sense. By the assertion (i) it follows that it is satisfied  pointwise  in $(\gamma ^{-1})\times (0, \infty)$.

The property (iv) is a direct consequence of property $(v)$ in Proposition \ref{S5PropU}. More precisely, since $s_0\in (0, \Re e(\sigma _2))$, it follows that $s<\Re e(\sigma _2+\sigma _1)\equiv 2$ and then, using the Lebesque's convergence Theorem and property $(v)$ of Proposition \ref{S5PropU}:
\bqns
\lim _{ \gamma t\to 1^+ }\omega (t, x)=\mathcal M _{ s_0 }^{-1}(U(\gamma ^{-1}))(x)=u(\gamma ^{-1}, x).
\eqns

In order to prove (v) we use the definition of $\omega $ and deformation of the contour integration. For $t>\gamma ^{-1}$ fixed and $x\to \infty$ we have:
\bqns
\omega (t, x)&=&
\frac{\gamma \Gamma \left(1+ \frac {\sigma _2} {\gamma }\right)}
{\Gamma \left( \frac {\sigma _1} {\gamma }\right)\Gamma \left(1+ \frac {\sigma _2-\sigma _1} {\gamma }\right)}\frac {(\gamma t-1)^{\frac {\sigma _1} {\gamma }-1}} {(\gamma t)^{\frac {\sigma _2} {\gamma }}}x^{-\sigma _2-\gamma }+
\frac {1} {2i\pi }\int  _{ \Re e s=s_* }U(t, s)x^{-s}ds
\eqns
for $s_*>\sigma _2+\gamma $. Then,
\bqns
\omega (t, x)&=&
\frac{\gamma \Gamma \left(1+ \frac {\sigma _2} {\gamma }\right)}
{\Gamma \left( \frac {\sigma _1} {\gamma }\right)\Gamma \left(1+ \frac {\sigma _2-\sigma _1} {\gamma }\right)}\frac {(\gamma t-1)^{\frac {\sigma _1} {\gamma }-1}} {(\gamma t)^{\frac {\sigma _2} {\gamma }}}x^{-\sigma _2-\gamma }+
\mathcal O\left(x^{-s_*} \right),\,\,\,x\to \infty
\eqns
and (v) follows.  
\qed
\vskip 0.15cm
\noindent

\begin{remark}
\label{S4Eexpl}
It is possible to obtain an explicit expression of $\omega (t, x)$ for all $\gamma t>1$ and $x>0$ by deforming the integration contour and the residue's Theorem.
For $\gamma x^\gamma t<1$ we must use the residues at $s=-m\gamma $, $m\in \N$:
\bqn
&&\omega (t, x)\equiv\frac {1} {2i\pi }\int  _{ \Re e s=s_0 }U(t, s)x^{-s}ds=
-\frac {\gamma\, (\gamma t)^{\frac {\sigma _1} {\gamma }}} {\Gamma \left( \frac {\sigma _1} {\gamma }\right)\Gamma \left(1- \frac {\sigma _1-\sigma _2} {\gamma }\right)}\times \nonumber\\
&&\times \sum_{ m=0 } ^\infty \frac {(\gamma x^\gamma t)^{m}\Gamma \left(1+m+\frac {\sigma _2 } {\gamma } \right)(-1)^{m+1}}
 {\Gamma (m+1)} F\left( 1-\frac {\sigma _1} {\gamma }, -m- \frac {\sigma _1} {\gamma }; 1-\frac {\sigma _1-\sigma _2} {\gamma }; \frac {1} {\gamma t}\right)\label{S5.5E2}
\eqn

If $\gamma x^\gamma t>1$ we use the poles at $s=\sigma _2+\gamma (m+1)$, $m\in \N$ and obtain
\bqn
&&\omega (t, x)\equiv\frac {1} {2i\pi }\int  _{ \Re e s=s_0 }U(t, s)x^{-s}ds=
-\frac {\gamma\, (\gamma t)^{\frac {\sigma _1-\sigma _2} {\gamma }-1}x^{-\sigma _2-\gamma }} {\Gamma \left( \frac {\sigma _1} {\gamma }\right)\Gamma \left(1- \frac {\sigma _1-\sigma _2} {\gamma }\right)}\times \label{S5.5E1}\\
&&\times \sum_{ m=0 } ^\infty \frac {(\gamma x^\gamma t)^{-m}\Gamma \left(1+m+\frac {\sigma _2 } {\gamma } \right)(-1)^{m+1} }
 {\Gamma (m+1)}F\left( 1-\frac {\sigma _1} {\gamma }, 1+m+ \frac {\sigma _2-\sigma _1} {\gamma }; 1-\frac {\sigma _1-\sigma _2} {\gamma }; \frac {1} {\gamma t}\right).\nonumber 
\eqn

The function $\omega $ may then be written as follows:
\bqns
\omega (t, x)=\frac{\gamma (\gamma t)^{\frac {\sigma _1} {\gamma }}H(t, \gamma tx^\gamma )}{\Gamma \left( \frac {\sigma _1} {\gamma }\right)\Gamma \left(1+ \frac {\sigma _2-\sigma _1} {\gamma }\right)},\,\,\,\forall t>\gamma ^{-1},\,\,\forall x>0,
\eqns
where
\bqns
H(t, z)=
\left\{
\begin{split}
&\sum_{ m=0 } ^\infty \frac {(-z)^{m}\Gamma \left(1+m+\frac {\sigma _2 } {\gamma } \right)}
 {\Gamma (m+1)}F\left( 1-\frac {\sigma _1} {\gamma }, -m- \frac {\sigma _1} {\gamma }; 1-\frac {\sigma _1-\sigma _2} {\gamma }; \frac {1} {\gamma t}\right),\,z\le 1\\
&z^{-\frac {\sigma _2} {\gamma }-1 }\sum_{ m=0 } ^\infty \frac {\Gamma \left(1+m+\frac {\sigma _2 } {\gamma } \right)}
 {(-z)^{m}\Gamma (m+1)} F\left( 1-\frac {\sigma _1} {\gamma }, 1+m+ \frac {\sigma _2-\sigma _1} {\gamma }; 1-\frac {\sigma _1-\sigma _2} {\gamma }; \frac {1} {\gamma t}\right),\,z\ge 1.
\end{split}
\right.
\eqns
\end{remark}

\noindent
\textbf{Proof of  Theorem \ref{S1MainTh3}.} We claim that the measure defined as
\bqn
\label{S1EGlobal}
w(t)=
\left\{
\begin{split}
&u(t),\,\,\,\hbox{if}\,\,\,t\in (0, \gamma ^{-1})\\
&\omega (t),\,\,\,\hbox{if}\,\,\,t\ge \gamma ^{-1}.
\end{split}
\right.
\eqn
satisfies all the requirements. 
In order to prove that $w$ is a global weak solution of  \eqref{eq:croisfrag2},\eqref{S1EK0} we must prove  that for all $\varphi \in C^1_c((0, \infty)\times (0, \infty))$:
\bqn
&&\int _0^\infty\int _0^\infty\left(\frac {\partial \varphi } {\partial t}+x^{\gamma +1}\frac {\partial \varphi } {\partial  x}+x^\gamma \varphi  (t, x) \right)w(t, x)dxdt=\nonumber\\
&&\hskip 1cm -\theta \int _0^\infty\int _0^\infty \varphi (t, x) \int _0^\infty w(t, y)y^{\gamma -1}dy\varphi (t, x)dxdt.\label{S6WeakS1}
\eqn
This is easily shown by  splitting the time  integrals  in  \eqref{S6WeakS1} in the two domains $(0, \gamma ^{-1})$ and $(\gamma ^{-1}, \infty)$, use \eqref{S3E565} in Remark \eqref{S3E564}, assertion (iii) of Proposition \ref{S6Propo2} and  the continuity of $w$ at $t=\gamma ^{-1}$.

Let us prove now   $w \in  \mathscr{C}((0, \infty); E' _{ 0, \Re e(\sigma _2)+\gamma  })$. By Theorem  \ref{S1MainTh1}, for all $t\in [0, \gamma ^{-1})$ the function $\mathcal M_{w}(t)=\mathcal M_{u}(t)$ is analytic  and bounded  in $\mathscr{S}(0, \infty)$. By Theorem \ref{S1MainTh2} , for $t>\gamma ^{-1}$, $\mathcal M_{w}(t)=\mathcal M_{\omega }(t)$ is analytic  and bounded  in $\mathscr{S}(0, \Re e(\sigma _2)+\gamma )$. But, for $t=\gamma ^{-1}$,  $\lim _{ t\to \gamma ^{-1} }\mathcal M_{w}(\gamma ^{-1}, s)$ is analytic and bounded only on
$\mathscr{S}(0, \nu)$. Then,  $\mathcal M_{w}(t)$ is analytic and bounded in $\mathscr{S}(0, \nu)$. Again by Theorem \ref {S1MainTh1} and Theorem \ref{S1MainTh2}, $\mathcal M_{w}\in \mathscr{C}((0, \infty)\times \mathscr{S}(0, \nu)$. Therefore $w \in  \mathscr{C}((0, \infty); E' _{ 0, \nu })$ using Theorem 11.10.1 in  \cite{ML}.

The properties \eqref{S1E198765} for $t>\gamma ^{-1}$ follow directly from Theorem \ref{S1MainTh2} since $w=\omega $ for $t>\gamma ^{-1}$. On the other since $u=w$ in $t\in (0, \gamma ^{-1})$ it follows from  \eqref{S3EPO1}-\eqref{S3EPO3} that $w$ also satisfies \eqref{S1E198765}  in that interval of time. \qed

\subsection{Uniqueness of the extension to $t\ge\gamma ^{-1}$.}
We are now concerned  with the question of uniqueness of global solutions to  \eqref{eq:croisfrag2},\eqref{S1EK0},\eqref{eq:croisfrag2data}. After the existence and uniqueness of a local solution $u$  on $(0, \gamma ^{-1})$ proved in Theorem \ref{S1MainTh1}, and since, by Corollary  \ref{S1Cor}, this local solution has a limit as $\gamma t\to 1^-$, this question is reduced in some sense to the uniqueness of the solutions of   \eqref{eq:croisfrag2},\eqref{S1EK0} with initial data $u\left( \gamma ^{-1}\right)$. How is this on the side of the Mellin variables?

By  general properties of hypergeometric functions,  the limit when $\gamma t\to 1^-$ of  $\mathcal M _{ u }$ only exists for $\Re e(s)<2$. Therefore, the  data  at $t=\gamma ^{-1}$ of  $\mathcal M _{ \omega  }$ is  only defined   for $\Re e(s)<2$. On the other hand, $\mathcal  M_{ \omega  }$  is meromorphic,  with a countable set of poles located at $s=-m\gamma$ and $s=\sigma _2+(m+1)\gamma$ for $m\in \N$. Since, in order to uniquely determine  $\mathcal  M_{ \omega  }$, we need the data  to be given in a strip of width strictly larger than $\gamma $, when $\gamma >2$ this forces to use an argument  of uniqueness in a strip where $\mathcal  M_{ \omega  }$ has  a pole. That is why  the hypothesis in Theorem \ref{S1MainTh2} are given in terms of $s\mathscr  P_{ \omega  }(t, s)$ on the strip $\mathscr{S}(-\gamma , \varepsilon)$, where $\mathscr P_{ \omega  }$ is the analytic extension of $\mathcal  M_{ \omega  }$ to $\mathscr{S}(-\gamma , \varepsilon)$.
\begin{theorem}
\label{S1MainTh2} 
For all $\gamma >0$ and $\theta >0$, the measure $\omega $ defined in \eqref{S6E642} is the unique real valued  weak solution $\omega \in  \mathscr{D}_1'((\gamma ^{-1}, \infty)\times (0, \infty))$ of  \eqref{eq:croisfrag2},\eqref{S1EK0} for all $t>\gamma ^{-1}$ such that its Mellin transform $\mathcal M_\omega$  satisfies the following properties: 
\bqn
\forall T>0, \exists \varepsilon >0:\hskip -0.5cm && \mathcal M_\omega(t)  \,\,\hbox{is analytic on} \,\,\mathscr{S}(0, \varepsilon),\,\,\forall t\in (\gamma ^{-1}, T) \nonumber\\
&&\forall t\in (0, T), \mathcal M_\omega(t, s)\,\, \hbox{has an extension}\,\,\mathscr{P} _{ \omega  },\,\hbox{to}\,\,(-\gamma , \varepsilon )\label{S1E2400}\\
&&\mathscr{P}_\omega\,\,\, \hbox{satisfies}\,\,  \eqref{S2eqconst} ,\,\,\hbox{on}\,\,(\gamma ^{-1}, T)\,\,\hbox{for some}\,\, s_*\in (-\gamma , \varepsilon)\label{S1E24007}\\
&&\forall t\in (0, T), s\mathscr{P} _{ \omega  }(t, s)\,\,\hbox{is analytic in}\,\,(-\gamma , \varepsilon )\\
&&s\mathscr{P} _{ \omega  }\in \mathscr{C}( [\gamma ^{-1},\infty)\times \mathscr{S}(-\gamma , \varepsilon))\label{S1E2401}\\ 
&&\sup\left\{|s\mathscr{P} _{ \omega  }(t, s)|;\,\, s\in \mathscr{S}(-\gamma  , \varepsilon),\,\,t\in \left[\gamma ^{-1}, T\right] \right\}<\infty \label{S1E2402}\\
&& \mathscr{P}_\omega(\gamma ^{-1}, s)=\Omega (\gamma ^{-1}, s),\,\,\,\forall s\in \mathscr{S}(-\gamma , \varepsilon).\label{S1E2403}
\eqn
Moreover this solution is such that:
\bqn
&&\omega \in C^\infty((\gamma ^{-1},\infty)\times (0, \infty))\cap  C([\gamma ^{-1},\infty)\times (0, \infty))\label{S1Ereg}.
\eqn
For all $t>\gamma ^{-1}$, $\omega $   satisfies  \eqref{eq:croisfrag2},\eqref{S1EK0} pointwise and its Mellin transform $\mathcal M_\omega $ is such that:
\bqn
&&\hskip -0.6cm \mathcal M_\omega  \,\,\hbox{is analytic  on}\,\,\, (\gamma ^{-1},\infty)\times \mathscr{S}(0, \Re e(\sigma _2)+\gamma ),\label{S4ES1014Bis}\\
&&\hskip -0.6cm \mathcal M_\omega (\gamma ^{-1},s)=
\Omega (\gamma ^{-1}, s)\,\,\, \forall s\in \C,\,\,\,\Re e(s)<2 \label{S1Ebw0+}
\\
&&\hskip -0.6cm \mathcal M_\omega \in \mathscr{C}( [\gamma ^{-1},\infty)\times \mathscr{S}(0, \nu)).\label{S4ES1015Bis}
\eqn
\end{theorem}

\begin{remark}
The condition \eqref{S1E2402} is not satisfied in general by $\mathcal M_u(t)$ for $\gamma t<1$, but  is satisfied  by $\mathcal M_u(\gamma ^{-1})$, the initial data of $\omega$ at $t=\gamma ^{-1}$ , as it follows using \eqref{S3EPO5} and Stirling's formula.
\end{remark}

\noindent
\textbf{Proof of  Theorem \ref{S1MainTh2}.} By  Proposition \ref{S5PropU}-(ii),
$\mathcal M _{ \omega  }(t,s)=U(t, s)$ for all $s\in \mathscr{S}(0, \Re e(\sigma _2)+\gamma )$. It follows from Proposition \ref{S5PropU} and Proposition \ref{S6Propo2} that $\omega $ satisfies 
\eqref{S1Ereg}-\eqref{S4ES1015Bis}. On the other hand, the function $\mathcal M _{ \omega  }(t,s)$ has a meromorphic  extension to the complex plane, given by  $U(t, s)$ that, by  Proposition \ref{S5PropU}-(ii) satisfies \eqref{S1E24007}. By  Proposition \ref{S5PropU}-(i), $U$ as a simple pole at $s=0$ and $sU(t, s)$ is analytic on $(-\gamma , \Re e(\sigma _2)+\gamma )$. By (i) and (iv) of that same Proposition, $sU(t, s)$ satisfies
\eqref{S1E2401} and it satisfies \eqref{S1E2402} and \eqref{S1E2403} by points (iii) and (iv). 

We prove  now the uniqueness of weak solutions satisfying  \eqref{S1E2400}--\eqref{S1E2403}. Suppose that two such solutions $\omega _1$ and $\omega _1$ exists and let $\mathscr{P}_{ \omega _1 }$, $\mathscr{P} _{ \omega _2 }$ be the extensions of their Mellin transforms. Then, the two functions 
$W_1(t, s)=\mathscr{P} _{ \omega _1 }(t-\gamma ^{-1}, s)$ and  $W_2(t, s)=\mathscr{P} _{ \omega _2 }(t-\gamma ^{-1}, s)$ satisfy the hypothesis of Theorem \ref{S4TU2} and are  then equal.
We deduce in particular that  $\mathcal M _{ \omega_1  }(t,s)=\mathcal M _{ \omega_2  }(t,s)$ for $t\in (\gamma ^{-1}, T)$ and $s\in \mathscr{S}(0, \varepsilon )$ which is a contradiction.

Suppose now that $\omega $ is complex valued. Since the coefficients of the equation  \eqref{eq:croisfrag2},\eqref{S1EK0}  are real, the conjugate $\overline \omega $ is also a solution, with the same initial data. Moreover, just by definition, its Mellin transform $\mathcal M _{ \overline \omega  }$ is such that  $\mathcal M _{ \overline \omega  }(t, s)=\overline{\mathcal M _{ \omega  }(t, \overline{s})}$ and $ s\mathcal M _{ \overline \omega  }(t, s)=\overline{\overline{s}\mathcal M _{ \omega  }(t, \overline{s})}$ for $s\in \mathscr{S}(0, \varepsilon )$. By hypothesis, for all $t>\gamma ^{-1}$, $\mathcal M _{ \omega  }(t, s)$ has an extension $\mathscr{P} _{ \omega  }(t, s)$ such that  the function $h(t, s)=s\mathscr{P} _{ \omega  }(t, s)$  is analytic in $s\in \mathscr{S}(-\gamma , \varepsilon)$. Therefore, $\mathscr{Q} _{\overline \omega  }(t, s)=\overline{\mathscr{P} _{ \omega  }(t, \overline{s})}=$ is an  extension of 
$\mathcal M _{ \overline \omega  }(t, s)$ such that  $\overline{h(t, \overline {s})}=\overline s\mathscr{Q} _{ \overline \omega  }(t, s)$ is analytic in the domain:
$$
\overline{\mathscr{S}}=\left\{s\in \C;\,\,\overline{s}\in \mathscr{S}(-\gamma , \varepsilon) \right\}\equiv \mathscr{S}(-\gamma , \varepsilon).
$$
Moreover, $\mathscr{Q} _{ \overline \omega  }$  satisfies the conditions \eqref{S1E2401}-\eqref{S1E2403} since, by hypothesis, so does  $\mathscr{P} _{ \omega  }$. We deduce that $u=\overline u$ by the uniqueness property that has been proved just above and $u$ is then real valued. This ends the proof of Theorem \ref{S1MainTh2}.
\qed
The following Proposition is used in the proof of Theorem \ref{S1MainTh2BIS},  
\begin{proposition}
\label{S1MainTh2Aux} Suppose that $\gamma (0, 2)$.
There exists a unique real valued  weak solution $\omega \in  \mathscr{D}_1'((\gamma ^{-1}, \infty)\times (0, \infty))$ of  \eqref{eq:croisfrag2},\eqref{S1EK0} for all $t>\gamma ^{-1}$ such that for some  $\rho\in (0, 2-\gamma )$ and all $T>\gamma ^{-1}$ its Mellin transform $M_\omega$ solves  \eqref{S2eqconst} and  satisfies:
\bqn
&&\forall t\in (\gamma ^{-1}, T),\mathcal M_\omega(t, s)\,\,\,  \hbox{is analytic on }\,\mathscr{S}(\rho-\delta , \rho +\gamma+\delta  ),\,\,\hbox{for some}\,\,\delta >0,\label{S1E2400T}\\
&&\mathcal M_\omega(t, s)\,\,\hbox{is continuous  on}\,\,\, [\gamma ^{-1},\infty)\times \mathscr{S}(\rho-\delta , \rho +\gamma+\delta  )\label{S1E2401T}\\ 
&&\sup\left\{|\mathcal M_\omega(s, t)|;\,\, s\in\mathscr{S}(\rho-\delta  , \rho +\gamma+\delta  ),\,\,t\in \left[\gamma ^{-1}, T\right] \right\}<\infty \label{S1E2402T}\\
&& \mathcal M_\omega(\gamma ^{-1}, s)=\Omega (\gamma ^{-1}, s),\,\,\,\forall s\in \mathscr{S}(\rho, \rho +\gamma ).\label{S1E2403T}
\eqn
This solution is the function $\omega $ obtained in Theorem \ref{S1MainTh2}.
 \end{proposition} 
\textbf{Proof of Proposition  \ref{S1MainTh2Aux}.} By Theorem \ref{S1MainTh2}  and the hypothesis $\gamma \in (0, 2)$, the function $\omega $ obtained in Theorem \ref{S1MainTh2}  satisfies \eqref {S1E2400T}-\eqref{S1E2403T}. Suppose that $\omega _1$ and $\omega _2$ are two different weak solution in
 $\mathscr{D}_1'((\gamma ^{-1}, \infty)\times (0, \infty))$ of  \eqref{eq:croisfrag2},\eqref{S1EK0} satisfying
\eqref {S1E2400T}-\eqref{S1E2403T}. Then the functions  $W_1(t, s)=\mathcal M _{\omega _1  }(t+\gamma ^{-1}, s)$ and $W_2(t, s)=\mathcal M _{\omega _2  }(t+\gamma ^{-1}, s)$ satisfy the hypothesis of  Theorem \ref{S4TU} for all $T>0$ for $W_0(s)=U(\gamma ^{-1}, s)$ and are then equal. This contradiction concludes the proof.
\qed
\vskip 0.15cm
We now prove Theorem \ref{S1MainTh2BIS} whose hypothesis are simpler than \eqref{S1E2400}-\eqref{S1E2403} in Theorem \ref{S1MainTh2}, but that requires the condition $\gamma \in (0, 2)$.
\vskip 0.15cm
\noindent
\textbf{Proof of Theorem \ref{S1MainTh2BIS}.} Suppose first that $T<\gamma^{-1}$. By Theorem \ref{S1MainTh1} the measure $u$  is a weak non negative solution of \eqref{eq:croisfrag2},\eqref{S1EK0} such that 
$u\in \mathscr{C}([0, \gamma ^{-1}); E' _{0, q })$ for all $q>0$ and satisfies \eqref{S1E2400Bis} .  If we suppose now that $u_1$ and $u_2$ are two different solutions satisfying these conditions, $\mathcal M _{ u_1 }$ and $\mathcal M _{ u_2 }$  would both satisfy the hypothesis of  Theorem \ref{S4TU} for all $T\in (0, \gamma ^{-1})$  and therefore would be equal. This contradiction proves the uniqueness for $T\in (0, \gamma ^{-1})$. 

Suppose now that $T>\gamma ^{-1}$ and $\tilde w$ satisfies the hypothesis of Theorem \ref{S1MainTh2BIS}.  We already know by the previous step that $\tilde w=u$ for $t\in (0, \gamma ^{-1})$. Let us prove that $\tilde w=\omega $ for $t\in (\gamma ^{-1}, T)$. By hypothesis $\tilde w \in  \mathscr{C}((0, T); E' _{ \rho -\delta , \rho +\gamma +\delta })$  for some $\rho\in (0, 2-\gamma )$, $\delta >0$ and is a weak solution of \eqref{eq:croisfrag2},\eqref{S1EK0}. Since 
$x_+^{s-1}\in E _{ p, q }$ for all $s\in \C$ such that $\Re e(s)\in (p, q)$ and $\mathscr{D} (0, \infty)$ is dense  in $E _{ p, q }$ for all $p<q$, we deduce that $\mathcal M _{\tilde  w  }$ solves the  equation \eqref{S2eqconst}. By \eqref{S1E2400Bis}, it also satisfies  \eqref{S1E2400T}-\eqref{S1E2402T}. Since, by the continuity of $\mathcal M _{ \tilde w }$ at $t=\gamma ^{-1}$ we also have $\mathcal M _{w  }(\gamma ^{-1}, s)=\mathcal M _{u  }(\gamma ^{-1}, s)=\Omega ((\gamma ^{-1}, s))$ for all $s\in \mathscr{S}(\rho , \rho +\gamma )$. , it follows that $\tilde w$ satisfies all the hypothesis of Proposition \ref{S1MainTh2Aux}.
We deduce from that Proposition that $\tilde w=\omega $ for $t\in (\gamma ^{-1}, T)$. We deduce by continuity that  $\tilde w=w$ for $t\in (0, T)$.
\qed

 \section{$\gamma >0$. Non existence of non negative solutions for large time.}
\label{sign}
The measure $w $ defined in Theorem  \ref{S1MainTh2} is a global weak solution of  \eqref{eq:croisfrag2},\eqref{S1EK0} for all values of the parameter $\theta$ and all possible values of the roots $\sigma _1$ and $\sigma _2$. However, as  we prove in this Section,  it is not always a non negative solution. By  uniqueness  of the possible extensions of the local solution $u$, as stated in Theorem  \ref{S1MainTh2}, it follows that it can not be extended to a suitable weak solution beyond  $t=\gamma ^{-1}$. \\
Our next result is concerned with the sign of the solution $\omega $ obtained in Theorem \ref{S1MainTh2}, and the possible extension of the local solution $u$ to a global non negative solution.

\begin{theorem}
\label{S1MainTh4} 
If $\theta>1$, the solution $w$ obtained in Theorem \ref{S1MainTh2}  is not always non negative for $t>\gamma ^{-1}$. 
\end{theorem}
\vskip 0.15cm
\noindent
\textbf{Proof of Theorem \ref{S1MainTh4}}
If we define the two following functions of $t>0$:
\bqns
A(t)&=&
\gamma \frac{ (\gamma t)^{\frac {\sigma _1-\sigma _2} {\gamma }-1}}{\Gamma \left( \frac {\sigma _1} {\gamma }\right)\Gamma \left(1+ \frac {\sigma _2-\sigma _1} {\gamma }\right)}\\
H(t)&=&\Gamma \left(1+\frac {\sigma _2 } {\gamma } \right) F\left( 1-\frac {\sigma _1} {\gamma }, 1+ \frac {\sigma _2-\sigma _1} {\gamma }; 1-\frac {\sigma _1-\sigma _2} {\gamma }; \frac {1} {\gamma t}\right)=\Gamma \left(1+\frac {\sigma _2 } {\gamma } \right)\left(1 - \frac {1} {\gamma t}\right)^{-1+\frac {\sigma _1} {\gamma }}
\eqns
 the right hand side of \eqref{S5.5E1} may  be written as follows:
\bqns
\omega (t, x)&=&A(t)x^{-\sigma _2-\gamma }\left(H(t)+B(t, x) \right)\\
B(t, x)&=&-(\gamma x^\gamma t)^{-1}\left(\sum_{ m=0 }^\infty \frac {(\gamma x^\gamma t)^{-m}\Gamma \left(2+m+\frac {\sigma _2} {\gamma } \right)(-1)^m} {\Gamma (2+m)}\times \right.\\
&&\hskip 5cm \times \left.F\left(1-\frac {\sigma _1} {\gamma }, 2+m+\frac {\sigma _2-\sigma _1} {\gamma }, \frac {1} {\gamma t} \right)  \right),
\eqns
where, by \eqref{S1EKS024}, $A(t)\ge 0$ and $H(t)\ge 0$ for all $t>\gamma ^{-1}$.
By  15.7.2 in \cite{AS}, there exists a constant $C=C(\sigma _1, \sigma _2, \gamma )$ such that for all $t>\frac {1} {\gamma }$ and $m\ge 0$:
$$
\left|F\left(1-\frac {\sigma _1} {\gamma }, 2+m+\frac {\sigma _2-\sigma _1} {\gamma }, \frac {1} {\gamma t} \right) \right|\le C e^{\left(2+m+\frac {\sigma _2-\sigma _1} {\gamma }\right)\frac {1} {\gamma t}}
$$
Then, for all $t>\gamma ^{-1}$, all $x>1$ and $m\ge 0$
\bqns
\left|\frac {(\gamma x^\gamma t)^{-m}\Gamma \left(2+m+\frac {\sigma _2} {\gamma } \right)(-1)^m} {\Gamma (2+m)}F\left(1-\frac {\sigma _1} {\gamma }, 2+m+\frac {\sigma _2-\sigma _1} {\gamma }, \frac {1} {\gamma t} \right) \right|\le
Ce^{-m\left( \log (\gamma x^\gamma t)-\frac {1} {\gamma t}\right)}\times \\
\times \frac {\Gamma \left(2+m+\frac {\sigma _2} {\gamma }\right)} {\Gamma (2+m)}.
\eqns
But, from  Stirling's formula:
$$
\frac {\Gamma \left(2+m+\frac {\sigma _2} {\gamma }\right)} {\Gamma (2+m)}=\mathcal O\left(\left( 2+m+\frac {\sigma _2} {\gamma }\right)^{\frac {\sigma _2} {\gamma }} \right),
\,\,\,m\to \infty.
$$
and therefore:
\bqns
\left|\sum_{ m=0 }^\infty \frac {(\gamma x^\gamma t)^{-m}\Gamma \left(2+m+\frac {\sigma _2} {\gamma } \right)(-1)^m} {\Gamma (2+m)}F\left(1-\frac {\sigma _1} {\gamma }, 2+m+\frac {\sigma _2-\sigma _1} {\gamma }, \frac {1} {\gamma t} \right)\right|\le\\
C \sum_{ m=0 }^\infty e^{-m\left( \log (\gamma x^\gamma t)-\frac {1} {\gamma t}\right)}m^{\frac {\sigma _2} {\gamma }}.
\eqns
The series in the right hand side defines a bounded function on the domain $t\ge \gamma ^{-1}$, $x\ge R _{ \gamma  }$ for $R _{ \gamma  }>0$ fixed large enough to have
$$
 \log (\gamma R _{ \gamma  }^\gamma t)-\frac {1} {\gamma t}\ge \delta >0,\,\,\,\forall t\ge\frac {1} {\gamma }.
$$
We then deduce that, for every $t>\gamma ^{-1}$:
\bqn
\label{S6E100}
\omega (t, x)=A(t)x^{-\sigma _2-\gamma }\left(H(t)+\mathcal O\left( \frac {1} {x^\gamma  t}\right) \right),\,\,x\to \infty.
\eqn
Suppose now that $\theta>1$, from where $\sigma _2\in \C\setminus \R$, and fix  any  $t_0>\gamma ^{-1}$.  There exists  $R=R(t_0)$ large enough such that:
\bqns
A(t_0)(H(t_0)+B(t, x))\ge  A(t_0)\left(H(t_0-|B(t, x)|\right)\\
\ge \frac {1} {2}A(t_0)H(t_0)>0,\,\,\,\forall x>R.
\eqns
Since $\sigma _2\in \C\setminus \R$, the function $x^{-\sigma _2-\gamma }$ is oscillatory.  Therefore by \eqref{S6E100}, $u(t_0, x)$ can not remain non negative for all $x>R$. 
\qed
\vskip 0.15cm
\begin{remark}
\label{S5Etheta1}
If $\theta=1$, then  $\sigma _1=\sigma _2=1$ and our final argument in the proof of Theorem  \ref{S1MainTh4} fails.
\end{remark}
\noindent
\textbf{Proof of Theorem \ref{S1MainTh4Bis}.} Suppose that there exists a non negative weak solution $\tilde w $ satisfying \eqref{S1E42},\eqref{S1E40} %and \eqref{S1E41}
 for some $T>\gamma ^{-1}$. Since $u\ge 0$, it follows from \eqref{S1E40}  that $u\in \mathscr{C}( [0, T); E' _{ \rho -\delta , \rho +\gamma +\delta  })$. The function, $\mathcal M _{ u }(t)$ is then well defined and analytic on $\mathscr{S}(\rho -\delta , \rho +\gamma +\delta )$ and satisfies \eqref{S2eqconst} for any $s_*\in (\rho , \rho +\gamma )$. By \eqref{S1E40} again, we deduce that $\mathcal M_u$ satisfies \eqref{S1E2400Bis}. Then  we have, by  Theorem \ref{S1MainTh2BIS},  $\tilde w(t)=w(t)$ for $t\in (T, \gamma ^{-1})$. But, by Theorem \ref{S1MainTh4}, $w$ is not non negative on $(0, \infty)$ when $t>\gamma ^{-1}$, and this contradiction concludes the proof. \qed

The following  Theorem on non existence of global solutions follows from the uniqueness result of Theorem \ref{S1MainTh2}  and Theorem \ref{S1MainTh4}, in the same way as Theorem \ref{S1MainTh4Bis}  follows from Theorem \ref{S1MainTh2BIS} and Theorem \ref{S1MainTh4}.
\begin{theorem}
\label{S5T621}
 Suppose that  $\theta>1$ and $T>\gamma ^{-1}$. Then there is no possible extension of the local solution $u$ to a non negative weak solution $\tilde w \in 
 \mathscr{D}_1'((0, T)\times (0, \infty))$ of  \eqref{eq:croisfrag2},\eqref{S1EK0} satisfying the initial condition \eqref{S1E42} and  such $\mathcal M _{ \tilde w }$ satisfies the conditions
 \eqref{S1E2400}--\eqref{S1E2403}.
\end{theorem}

When $\theta \in (0, 1)$ the condition  \eqref{S1prop1B}  is satisfied and by Corollary 4.2 in \cite{BW} a global nonnegative solution $\mu $ exists. Although we do not know if $\mu =w$, the solution obtained in Theorem  \ref{S1MainTh3}, in general, we have:

\begin{proposition}
\label{finalrem}
Suppose that $\theta\in (0, 1)$ and $0<\gamma <\sqrt {1-\theta}$. Let $\mu $ be the global non negative solution of \eqref{eq:croisfrag2},\eqref{S1EK0} obtained in \cite{BW}  and $w$ the solution obtained in Theorem  \ref{S1MainTh3}. Then $\mu =w$.
\end{proposition}
\textbf{Proof.}
When $\theta \in (0, 1)$ condition  \eqref{S1prop1B}  is satisfied and  by Corollary 4.2 in \cite{BW} a global nonnegative solution $\mu $ exists. Moreover, it  follows from Lemma 4.3-i) and Corollary 4.2 in \cite{BW} that  $\mu  \in  \mathscr{C}([0, T); E' _{1,  \sigma _2})$  and  $\mathcal M _{ \mu }\in $ satisfies  \eqref{S1E2400Bis} for any $\rho \in (1, \sigma _2-\gamma )$.
By the uniqueness of such solutions proved in Theorem \ref{S1MainTh2BIS} it follows that $\mu=w$. \qed

\begin{remark}
The results in \cite{BW} are proved  for  general dislocation measures, for which the corresponding function $\Phi(s)$ could be  defined only for $\Re e s\ge 2$. When $k_0$ as in \eqref{S1EK0},  $\Phi (s)$ is  defined for all   $\Re e s>-1$ and the solutions in \cite{BW} may then be expected to  have moments of order $s$ in a larger interval  than  $\Re e(s)\in (1, \sigma _2)$. That could make Proposition \ref{finalrem} to be true under a weaker  condition than  $\gamma <\sqrt {1-\theta}$.
\end{remark}

\section{The case $\gamma <0$.}
\label{SGammaNegatif}
When $\gamma <0$ the function obtained in Section \ref{Explicit}
\bqns
\label{S2solmellin}
\Omega (t, s)=F\left(\frac {s-\sigma _1} {\gamma }, \frac {s-\sigma _2} {\gamma }, \frac {s} {\gamma }, \gamma t\right)
\equiv  (1-\gamma t)^{\frac {2-s} {\gamma }}F\left(\frac {\sigma _1} {\gamma }, \frac {\sigma _2} {\gamma  }, \frac {s} {\gamma }, \gamma t \right).
\eqns
is  a solution  of \eqref{S2eqconst} for all $t>0$. If,  to obtain a solution to \eqref{eq:croisfrag2},\eqref{S1EK0}, our purpose was still to take its inverse Mellin transform,  the inverse Mellin transform should  be defined along a vertical  integration curve contained in the half plane $\Re e(s)>0$, as in the case $\gamma >0$. 
But now the poles of $\Omega (t)$, namely  $s=-m\gamma, m\in \N$, are  non negative real numbers. It follows that  the  moments $\mathcal  M_u(t, r)$ of any solution $u$ of \eqref{eq:croisfrag2},\eqref{S1EK0} whose Mellin transform is $\Omega (t, s)$ will be bounded only  for $r$ in an interval  $(-m \gamma, -(m+1)\gamma  )$  for some $m\in \N$. For that reason we look for another solution of \eqref{S2eqconst}-\eqref{S2initial}

\subsection{Still another solution of \eqref{S2eqconst}-\eqref{S2initial}.}
The function 
\bqn
\label{S5EqV1}
V_2(s)=\frac {\gamma ^{\frac {s} {\gamma }}\Gamma \left(1-\frac {s} {\gamma  } \right)} 
{\Gamma \left(1-\frac {s-\sigma _1} {\sigma } \right)\Gamma \left(1-\frac {s-\sigma _2} {\sigma }\right) }
\eqn
satisfies:
\bqn
\label{S5EqV}
V_2(s+\gamma )=-\frac {(s-\sigma _1)(s-\sigma _2)} {s}V_2(s),\,\,\,\forall s\in \C\setminus \{s\in \C; s=-m\gamma ,\, m\in \N\}.
\eqn

For $\sigma _0>0$ fixed let $\widetilde{\mathscr{C} _{ \theta }}$ be the following curve in the complex plane:
\bqns
\widetilde{\mathscr{C} _{ \theta }}&=&\widetilde{\mathscr{C}_{1,\theta}}\cup \widetilde{\mathscr{C}_{2}}\\
\widetilde{\mathscr{C} _{ 1, \theta }}&=& \left\{s=\xi +i\zeta \in \C; \xi =\sigma _0+\theta \zeta ,\,\,\zeta \ge 0 \right\}\\
\widetilde{\mathscr{C}_{2}}&=& \left\{s=\xi +i\zeta \in \C; \xi =\sigma _0,\,\,\zeta <0 \right\}
\eqns
and define
\bqn
\label{S610}
U_2(t, s)=\frac {1} { \gamma V_2(s)}\int _{\widetilde{\mathscr{C} _{ \theta }}}\frac {(-t)^{\frac {\sigma -s} {\gamma }}V_2(\sigma )d\sigma } 
{\Gamma \left( 1+\frac {\sigma -s} {\gamma }\right)\left(e^{-\frac {2i\pi } {\gamma }(\sigma -s)}-1\right)}.
\eqn

\begin{proposition} 
\label{S6P1} For all $t\in (0, (-\gamma )^{-1})$, the integral in the right hand side of \eqref{S610} is absolutely convergent and defines the following meromorphic  function  $U_2(t)$ in the complex plane:
\bqn
U_2(t, s)&=&\Omega _1(t, s)-\Omega _2(t, s)\label{S6E20}\\
\Omega _1(t, s)&=&\frac {(-\gamma t)^{-\frac {s} {\gamma }}} {\left(e^{\frac {2i\pi } {\gamma }s}-1\right)}
 \frac {\gamma t\, \Gamma \left(1-\frac {s-\sigma _1 } {\gamma  }\right) \Gamma \left(1-\frac {s-\sigma _2 } {\gamma  }\right)} {\Gamma \left(\frac {\sigma _1 } {\gamma  }\right)\Gamma \left(\frac {\sigma _2 } {\gamma  }\right)\Gamma \left(1-\frac {s  } {\gamma  }\right)}
 \frac {F\left( 1-\frac {\sigma _1} {\gamma }, 1-\frac {\sigma _2} {\gamma }, 2-\frac {s} {\gamma }, \gamma t\right)} {\Gamma \left(2-\frac {s  } {\gamma  }\right)}\label{S6E21}\\
\Omega _2(t, s)&=&\frac {i} { 2\pi}F\left( \frac {s-\sigma _1} {\gamma }, \frac {s-\sigma _2} {\gamma }, \frac {s} {\gamma }, \gamma t\right)
=\frac {i} { 2\pi}(1-\gamma t)^{\frac {2-s} {\gamma }}F\left( \frac {\sigma _1} {\gamma }, \frac {\sigma _2} {\gamma }, \frac {s} {\gamma }, \gamma t\right)\label{S6E22}
\eqn
and
\bqn
\label{S6E201}
U_2\in \mathscr C([0, -\gamma ^{-1})\times \mathscr{S}(1-\gamma , \infty)).
\eqn

\end{proposition}
\textbf{Proof.} Using Stirling's formulas in the same way as in Section \ref{S6}, we obtain by  straightforward calculation that for all $t\in (0, (-\gamma )^{-1})$ and $s\in \C$, there exists a positive constant $C=C(s, \gamma )$:
\bqn
\label{S6E12}
\left|\frac {(-t)^{\frac {\sigma} {\gamma }}V_2(\sigma )d\sigma } 
{\Gamma \left( 1+\frac {\sigma -s} {\gamma }\right)\left(e^{-\frac {2i\pi } {\gamma }(\sigma -s)}-1\right)} \right|
\le 
\left\{
\begin{split}
&Ce^{\frac {\xi  } {\gamma }\log(-\gamma t)}e^{-\frac {\zeta  \pi } {\gamma }},\,\,\,v>0\\
&C e^{\frac {\xi } {\gamma }\log(-\gamma t)}e^{\frac {2\zeta \pi } {\gamma }},\,\,\,v<0.
\end{split}
\right.
\eqn
The right hand side is exponentially decaying as $\zeta \to -\infty$ for any $\xi $ fixed. On the other hand, suppose that $-\gamma t\in (0, \tau )$ with $\tau <1$. Then 
$\log (-\gamma t)<\log \tau <0$ and 
\bqns
0<\frac {\log \tau } {\gamma }<\frac {\log(-\gamma t)} {\gamma }.
\eqns
Therefore, if $\sigma =\xi +i\zeta \in \widetilde{\mathscr{C} _{ 1, \theta }}$, $\xi =\sigma _0+\theta \zeta $ and:
\bqns
\frac {\xi } {\gamma }\log(-\gamma t)-\frac {\zeta \pi } {\gamma }&=&\frac {(\sigma _0+\theta \zeta )} {\gamma }\log(-\gamma t)-\frac {\zeta \pi } {\gamma }=
\frac {\sigma _0} {\gamma }\log(-\gamma t)+\theta \zeta \frac {\log(-\gamma t)} {\gamma }-\frac {\zeta \pi } {\gamma }\\
&=&\frac {\sigma _0} {\gamma }\log(-\gamma t)+\zeta \left( \theta \frac {\log(-\gamma t)} {\gamma }-\frac {\pi } {\gamma }\right).
\eqns
If we choose now $\theta<0$ we deduce:
\bqns
\frac {\xi } {\gamma }\log(-\gamma t)-\frac {\zeta \pi } {\gamma }<\frac {\sigma _0} {\gamma }\log(-\gamma t)+\zeta \left( \theta\frac {\log \tau } {\gamma }-\frac {\pi } {\gamma }\right)
\eqns
and, if
\bqns
\theta <\frac {\pi } {\log \tau }
\eqns
the right hand side of \eqref{S6E12} is also exponentially decreasing as $\zeta \to \infty$ and $\sigma \in \mathscr{C _{ \theta }}$. The function under the integral  in \eqref{S610} is then absolutely integrable.

The integral may now be obtained using the method of residues.
The poles of the function to integrate are:
\bqns
\sigma =\gamma (m+1),\,\,\, \sigma =s+\gamma m,\,\,\,m\in \N.
\eqns
and then, for $-\gamma t\in (0, 1)$:
\bqns
U(t, s)=\frac {(-t)^{-\frac {s} {\gamma }}} {\gamma V(s)}\sum_{ m=0 } ^\infty\frac {(-t)^{m+1}\gamma ^{m+1} Res\left(\Gamma \left(1-\frac {\sigma } {\gamma  } \right); \sigma =\gamma (m+1)\right)}
 {\Gamma \left(\frac {\sigma _1} {\gamma  } -m\right)\Gamma \left(\frac {\sigma _2} {\sigma }-m\right)\Gamma \left( 2+m-\frac {s} {\gamma }\right)\left(e^{\frac {2i\pi } {\gamma }s}-1\right)}+\\
 +\frac {1} {\gamma V(s)}\sum_{ m=0 } ^\infty
 \frac {(-t)^{m}\gamma ^{\frac {s} {\gamma }+m} \Gamma \left(1-\frac {s } {\gamma  }-m\right) Res\left( \left(e^{\frac {2i\pi } {\gamma }{(\sigma -s)}}-1\right)^{-1}; \sigma =s+\gamma m\right)}
 {\Gamma \left(1-\frac {s-\sigma _1} {\gamma  } -m\right)\Gamma \left(1-\frac {s-\sigma _2} {\sigma }-m\right)\Gamma \left( 1+m\right)}.
\eqns 
We use that, for all $m\in \N$:
\bqns
Res\left(\Gamma \left(1-\frac {\sigma } {\gamma  } \right); \sigma =\gamma (m+1)\right)=\frac {(-1)^{m+1}\gamma } {\Gamma (m+1)}\\
Res\left( \left(e^{\frac {2i\pi } {\gamma }(\sigma -s)}-1\right)^{-1}; \sigma =s+\gamma m\right)=-\frac {i\gamma } {2\pi }
\eqns
and obtain, for $-\gamma t\in (0, 1)$:
\bqns
U(t, s)=\frac {(-t)^{-\frac {s} {\gamma }}} { V(s)\left(e^{\frac {2i\pi } {\gamma }s}-1\right)}\sum_{ m=0 } ^\infty\frac {(t)^{m+1}\gamma ^{m+1} }
 {\Gamma \left(\frac {\sigma _1} {\gamma  } -m\right)\Gamma \left(\frac {\sigma _2} {\sigma }-m\right)\Gamma \left( 2+m-\frac {s} {\gamma }\right)\Gamma (m+1)}-\\
-\frac {i\gamma ^{\frac {s} {\gamma }} } { 2\pi V(s)}\sum_{ m=0 } ^\infty
 \frac {(-t)^{m}\gamma ^{m} \Gamma \left(1-\frac {s } {\gamma  }-m\right)}
 {\Gamma \left(1-\frac {s-\sigma _1} {\gamma  } -m\right)\Gamma \left(1-\frac {s-\sigma _2} {\sigma }-m\right)\Gamma \left( 1+m\right)}.
\eqns 
The two series may be summed when  $-\gamma t\in (0, 1)$:
\bqns
\sum_{ m=0 } ^\infty\frac {(t)^{m+1}\gamma ^{m+1} }
 {\Gamma \left(\frac {\sigma _1} {\gamma  } -m\right)\Gamma \left(\frac {\sigma _2} {\sigma }-m\right)\Gamma \left( 2+m-\frac {s} {\gamma }\right)\Gamma (m+1)}=
 \frac {\gamma t\, F\left( 1-\frac {\sigma _1} {\gamma }, 1-\frac {\sigma _2} {\gamma }, 2-\frac {s} {\gamma }, \gamma t\right)} {\Gamma \left(\frac {\sigma _1 } {\gamma  }\right)\Gamma \left(\frac {\sigma _2 } {\gamma  }\right)\Gamma \left(2-\frac {s  } {\gamma  }\right)}\\
 \sum_{ m=0 } ^\infty
 \frac {(-t)^{m}\gamma ^{m} \Gamma \left(1-\frac {s } {\gamma  }-m\right)}
 {\Gamma \left(1-\frac {s-\sigma _1} {\gamma  } -m\right)\Gamma \left(1-\frac {s-\sigma _2} {\sigma }-m\right)\Gamma \left( 1+m\right)}
 =\frac {\Gamma \left(1-\frac {s } {\gamma  }\right)F\left( \frac {s-\sigma _1} {\gamma }, \frac {s-\sigma _2} {\gamma }, \frac {s} {\gamma }, \gamma t\right)
 } { \Gamma \left(1-\frac {s-\sigma _1 } {\gamma  }\right) \Gamma \left(1-\frac {s-\sigma _2 } {\gamma  }\right)}.
\eqns
We deduce:
\bqns
U_2(t, s)=\frac {(-t)^{-\frac {s} {\gamma }}} { V_2(s)\left(e^{\frac {2i\pi } {\gamma }s}-1\right)} \frac {\gamma t\, F\left( 1-\frac {\sigma _1} {\gamma }, 1-\frac {\sigma _2} {\gamma }, 2-\frac {s} {\gamma }, \gamma t\right)} {\Gamma \left(\frac {\sigma _1 } {\gamma  }\right)\Gamma \left(\frac {\sigma _2 } {\gamma  }\right)\Gamma \left(2-\frac {s  } {\gamma  }\right)}-\\
-\frac {i\gamma ^{\frac {s} {\gamma }} } { 2\pi V_2(s)}\frac {\Gamma \left(1-\frac {s } {\gamma  }\right)F\left( \frac {s-\sigma _1} {\gamma }, \frac {s-\sigma _2} {\gamma }, \frac {s} {\gamma }, \gamma t\right)
 } { \Gamma \left(1-\frac {s-\sigma _1 } {\gamma  }\right) \Gamma \left(1-\frac {s-\sigma _2 } {\gamma  }\right)}.
\eqns
and the explicit expression of $U_2(t)$ in  \eqref{S6E20}-\eqref{S6E22}  follows using the expression \eqref{S5EqV1} of $V_2(s)$. The property \eqref{S6E201} follows now from 
\eqref{S6E20}-\eqref{S6E22} and the well known properties of the Gamma and hypergeometric functions
\qed 

\subsection{The solution $v$.}
We wish now to define a solution $u$ of \eqref{eq:croisfrag2},\eqref{S1EK0} by means of a suitable inverse Mellin transform of $U_2(t)$. We have already obtained the inverse Mellin transform of $\Omega _2$ in Section \ref{S4}. Some useful properties of $\Omega _1$ are  now given in the next  Proposition.

\begin{proposition}
\label{S6P10}
For all $t\in (0, -\gamma ^{-1})$, the function $\Omega _1(t)$ is analytic for $s\in\C$ such that $\Re e(s)>1+\gamma $. For all $t\in (0, -\gamma^{-1})$ there exists a positive constant $C$ such that
\bqn
\label{S6E210}
\left|\Omega _1(t, s) \right|
\le 
\left\{
\begin{split}
&C(t)\,e^{\frac {2\pi \Im m(s) } {\gamma }},\,\,\,\Im m(s)>0\\
&C(t)\,(1+|\Im m (s)|)^{-1+\frac {2} {\gamma }},\,\,\,\Im m(s)<0.
\end{split}
\right.
\eqn
\end{proposition}
\textbf{Proof.} The estimate follows from the expression of $\Omega _1$ and Stirling's formula. \qed
\vskip 0.15cm

We deduce from the  Proposition \eqref{S6P10}  that we may set
\bqn
\label{S6sol1}
v(t, x)=\frac {1} {2i\pi }\int  _{ \Re e(s)=s_0 }x^{-s}U_2(t, s)ds,\,\,\,s_0>1+\gamma
\eqn
and, by classical results on the Mellin and inverse Mellin transform this expression defines a measure 
$v\in   \mathscr{C}\left([0, -\gamma ^{-1}); E' _{ 1+\gamma , \infty}\right)$.

\begin{proposition}
For all $t\in (0, -\gamma ^{-1})$, $U_2(t)$ is meromorphic on $\C$ with poles  located at:
\bqns
s=\sigma  _{ \ell }+(m+1)\gamma ,\,\,\ell=1, 2,\,\,m\in \N.
\eqns
The residues of $U_2(t)$ at these points are:
\bqns
Res\left(U(t, s); s=\sigma _{ 1}+\gamma (m+1)\right)&=&\frac {(-\gamma t)^{-\frac {\sigma _1} {\gamma }-(m+1)}} {\left(e^{\frac {2i\pi } {\gamma }\sigma _1}-1\right)}
 \frac {\gamma t\, \Gamma \left(-\frac {\sigma _1-\sigma _2 } {\gamma  }-m\right)(-1)^{m+1}}
  {\Gamma \left(\frac {\sigma _1 } {\gamma  }\right)\Gamma \left(\frac {\sigma _2 } {\gamma  }\right)\Gamma \left(-\frac {\sigma _1  } {\gamma  }-m\right)\Gamma (m+1)}\times \\
&&  \times 
 \frac {F\left( 1-\frac {\sigma _1} {\gamma }, 1-\frac {\sigma _2} {\gamma }, 1-\frac {\sigma _1} {\gamma }-m, \gamma t\right)} {\Gamma \left(1-\frac {\sigma _1 } {\gamma  }-m\right)}
 =:A_m(t)\\
 Res\left(U(t, s); s=\sigma _{ 2}+\gamma (m+1)\right)&=&\frac {(-\gamma t)^{-\frac {\sigma _2} {\gamma }-(m+1)}} {\left(e^{\frac {2i\pi } {\gamma }\sigma _2}-1\right)}
 \frac {\gamma t\, \Gamma \left(-\frac {\sigma _2-\sigma _1 } {\gamma  }-m\right)(-1)^{m+1}}
  {\Gamma \left(\frac {\sigma _1 } {\gamma  }\right)\Gamma \left(\frac {\sigma _2 } {\gamma  }\right)\Gamma \left(-\frac {\sigma _2  } {\gamma  }-m\right)\Gamma (m+1)}\times \\
&&\times  \frac {F\left( 1-\frac {\sigma _1} {\gamma }, 1-\frac {\sigma _2} {\gamma }, 1-\frac {\sigma _2} {\gamma }-m, \gamma t\right)} {\Gamma \left(1-\frac {\sigma _2  } {\gamma  }-m\right)} =:B_m(t).
\eqns
Moreover, $U_2$ satisfies   \eqref{S2eqconst} for $t>\gamma ^{-1}$ and $s$ such that $\Re e(s)>1+\gamma$.
\end{proposition}
\textbf{Proof.} For each $t\in (0, -\gamma ^{-1})$ the functions $\Omega _1(t)$ and $\Omega _2(t)$ are meromorphic on $\C$ with poles located respectively at
$s=-\gamma m,\,\,s=\sigma _{ \ell }+\gamma (m+1)$ and $s=-\gamma m$, $m\in \N$.
But, at poles $s=-m\gamma $ we have: 
\bqns
&&Res\left(\Omega _1(t), s=-\gamma m \right)= Res\left( \left(e^{\frac {2i\pi } {\gamma }{s}}-1\right)^{-1}; s =-\gamma m\right)\times \\
&&\times 
(-\gamma t)^{m}
 \frac {\gamma t\, \Gamma \left(1+\frac {\sigma _1 } {\gamma  }+m\right) \Gamma \left(1+\frac {\sigma _2 } {\gamma  }+m\right)} {\Gamma \left(\frac {\sigma _1 } {\gamma  }\right)\Gamma \left(\frac {\sigma _2 } {\gamma  }\right)\Gamma \left(1+m\right)}
 \frac {F\left( 1-\frac {\sigma _1} {\gamma }, 1-\frac {\sigma _2} {\gamma }, 2+m, \gamma t\right)} {\Gamma \left(2+m\right)}\\
&& =-\frac {i\gamma } {2\pi }
(-\gamma t)^{m}
 \frac {\gamma t\, \Gamma \left(1+\frac {\sigma _1 } {\gamma  }+m\right) \Gamma \left(1+\frac {\sigma _2 } {\gamma  }+m\right)} {\Gamma \left(\frac {\sigma _1 } {\gamma  }\right)\Gamma \left(\frac {\sigma _2 } {\gamma  }\right)\Gamma \left(1+m\right)}
 \frac {F\left( 1-\frac {\sigma _1} {\gamma }, 1-\frac {\sigma _2} {\gamma }, 2+m, \gamma t\right)} {\Gamma \left(2+m\right)}.
\eqns
On the other hand, since
\bqns
F\left( \frac {\sigma _1} {\gamma }, \frac {\sigma _2} {\gamma }, \frac {s} {\gamma }, \gamma t\right)=
\sum_{ n=0 } ^\infty \frac {\Gamma \left(\frac {\sigma _1 } {\gamma  }+n\right)\Gamma \left(\frac {\sigma _2 } {\gamma  }+n\right)\Gamma \left(\frac {s} {\gamma  }\right)(\gamma t)^n} 
{\Gamma \left(\frac {\sigma _1 } {\gamma  }\right)\Gamma \left(\frac {\sigma _2 } {\gamma  }\right)\Gamma \left(\frac {s} {\gamma  }+n\right)\Gamma (n+1)}.
\eqns
we deduce:
\bqns
Res\left(F\left( \frac {\sigma _1} {\gamma }, \frac {\sigma _2} {\gamma }, \frac {s} {\gamma }, \gamma t\right), s=-\gamma m \right)
=\gamma \sum _{n= m+1 }^\infty
 \frac {\Gamma \left(\frac {\sigma _1 } {\gamma  }+n\right)\Gamma \left(\frac {\sigma _2 } {\gamma  }+n\right)(\gamma t)^n(-1)^{m}} 
{\Gamma \left(\frac {\sigma _1 } {\gamma  }\right)\Gamma \left(\frac {\sigma _2 } {\gamma  }\right)\Gamma \left(-m+n\right)\Gamma (n+1)\Gamma (m+1)}.
\eqns
This series may still be summed,
\bqns
Res\left(F\left( \frac {\sigma _1} {\gamma }, \frac {\sigma _2} {\gamma }, \frac {s} {\gamma }, \gamma t\right), s=-\gamma m \right)
=\gamma \frac {(\gamma t)^{m+1}(-1)^m\Gamma \left(1+m+\frac {\sigma _1 } {\gamma  }\right)\Gamma \left(1+m+\frac {\sigma _2 } {\gamma  }\right)} {\Gamma \left(\frac {\sigma _1 } {\gamma  }\right)\Gamma \left(\frac {\sigma _2 } {\gamma  }\right)\Gamma (m+2)\Gamma (m+1)}\times \\
\times  F\left(1+m+ \frac {\sigma _1} {\gamma }, 1+m+ \frac {\sigma _2} {\gamma }, 2+m, \gamma t\right).
\eqns
We use now:
\bqns
F\left(1+m+ \frac {\sigma _1} {\gamma }, 1+m+ \frac {\sigma _2} {\gamma }, 2+m, \gamma t\right)=
(1-\gamma t)^{-m-\frac {2} {\gamma }}F\left(1- \frac {\sigma _1} {\gamma }, 1- \frac {\sigma _2} {\gamma }, 2+m, \gamma t\right)
\eqns
from where:
\bqns
Res\left(F\left( \frac {\sigma _1} {\gamma }, \frac {\sigma _2} {\gamma }, \frac {s} {\gamma }, \gamma t\right), s=-\gamma m \right)
=\gamma \frac {(\gamma t)^{m+1}(-1)^m\Gamma \left(1+m+\frac {\sigma _1 } {\gamma  }\right)\Gamma \left(1+m+\frac {\sigma _2 } {\gamma  }\right)} {\Gamma \left(\frac {\sigma _1 } {\gamma  }\right)\Gamma \left(\frac {\sigma _2 } {\gamma  }\right)\Gamma (m+2)\Gamma (m+1)}\times \\
\times  (1-\gamma t)^{-m-\frac {2} {\gamma }}F\left(1- \frac {\sigma _1} {\gamma }, 1- \frac {\sigma _2} {\gamma }, 2+m, \gamma t\right)
\eqns
and then,
\bqns
Res\left(\Omega _2(t), s=-\gamma m \right)&=&\frac {i} {2\pi }\gamma \frac {(\gamma t)^{m+1}(-1)^m\Gamma \left(1+m+\frac {\sigma _1 } {\gamma  }\right)\Gamma \left(1+m+\frac {\sigma _2 } {\gamma  }\right)} {\Gamma \left(\frac {\sigma _1 } {\gamma  }\right)\Gamma \left(\frac {\sigma _2 } {\gamma  }\right)\Gamma (m+2)\Gamma (m+1)}\times \\
&&\times F\left(1- \frac {\sigma _1} {\gamma }, 1- \frac {\sigma _2} {\gamma }, 2+m, \gamma t\right).
\eqns
Therefore, at $s=-m\gamma $, the residues of $\Omega _1(t)$ and $\Omega _2(t)$ are equal and therefore,  they cancel when combined to obtain the residue of $U(t)$.
On the other hand, the residues of $\Omega _1(t)$ at $s=\sigma _{ 2 }+\gamma (m+1)$:
\bqns
Res\left(\Omega _1(t, s); s=\sigma _{ 2 }+\gamma (m+1)\right)&=&\frac {(-\gamma t)^{-\frac {\sigma _2} {\gamma }-(m+1)}} {\left(e^{\frac {2i\pi } {\gamma }\sigma _2}-1\right)}
 \frac {\gamma t\, \Gamma \left(-\frac {\sigma _2-\sigma _1 } {\gamma  }-m\right)(-1)^{m+1}}
  {\Gamma \left(\frac {\sigma _1 } {\gamma  }\right)\Gamma \left(\frac {\sigma _2 } {\gamma  }\right)\Gamma \left(-\frac {\sigma _2  } {\gamma  }-m\right)\Gamma (m+1)}\times \\
&&  \times 
 \frac {F\left( 1-\frac {\sigma _1} {\gamma }, 1-\frac {\sigma _2} {\gamma }, 1-\frac {\sigma _2} {\gamma }-m, \gamma t\right)} {\Gamma \left(1-\frac {\sigma _2  } {\gamma  }-m\right)}
\eqns
and similarly:
\bqns
Res\left(\Omega _1(t, s); s=\sigma _{ 1}+\gamma (m+1)\right)&=&\frac {(-\gamma t)^{-\frac {\sigma _1} {\gamma }-(m+1)}} {\left(e^{\frac {2i\pi } {\gamma }\sigma _1}-1\right)}
 \frac {\gamma t\, \Gamma \left(-\frac {\sigma _1-\sigma _2 } {\gamma  }-m\right)(-1)^{m+1}}
  {\Gamma \left(\frac {\sigma _1 } {\gamma  }\right)\Gamma \left(\frac {\sigma _2 } {\gamma  }\right)\Gamma \left(-\frac {\sigma _1  } {\gamma  }-m\right)\Gamma (m+1)}\times \\
&&  \times 
 \frac {F\left( 1-\frac {\sigma _1} {\gamma }, 1-\frac {\sigma _2} {\gamma }, 1-\frac {\sigma _1} {\gamma }-m, \gamma t\right)} {\Gamma \left(1-\frac {\sigma _1 } {\gamma  }-m\right)}.
\eqns
Arguing as in the proof of Proposition \ref{S5PropU},  we deduce that $U_2$ satisfies   \eqref{S2eqconst} for $t>\gamma ^{-1}$ and $s$ such that $\Re e(s)>1+\gamma $.
\qed

\begin{theorem}
\label{S6T12345}
The measure  $u\in   \mathscr{C}\left([0, -\gamma ^{-1}); E' _{ 1+\gamma , \infty}\right)$ defined in \eqref{S6sol1} is a weak solution of \eqref{eq:croisfrag2},\eqref{S1EK0} 
on $t\in (0, -\gamma ^{-1})$. It satisfies
\bqn
\label{S6E500}
\int  _{ 0 }^\infty v(t, x)x^{s-1}dx=U_2 (t, s),\,\,\,\forall s \in \C;\,\,\Re e(s)>1+\gamma.
\eqn
If $\theta >1$, for all $t\in (0, -\gamma ^{-1})$ the measure $v(t)$ takes positive and negative values on $(0, \infty)$.
\end{theorem}
\textbf{Proof.} The identity \eqref{S6E500} follows from Theorem  11.10.1 in \cite{ML} and by classical properties of the Mellin transform, $u$ is a weak solution of 
\eqref{eq:croisfrag2},\eqref{S1EK0}.

It follows from \eqref{S6sol1} and the properties of $U_2(t)$ that as $x\to 0$:
\bqns
v(t, x)=\Re e\left(A_0(t)x^{-\sigma _2-\gamma }\right)+ \Re e\left(B_0(t)x^{-\sigma _1-\gamma }\right)+o\left(x^{-1-\gamma } \right),\,\,\,x\to 0
\eqns
If $\theta>1$, $\sigma _2=1+i\zeta $, $\sigma _1=1-i\zeta $ with $\zeta =\sqrt{\theta-1}$,
\bqns
x^{-\sigma _2-\gamma  }=x^{-1-i\zeta -\gamma }=x^{-1-\gamma }\left(\cos(\zeta \log x)-i\sin(\zeta \log x) \right)\\
x^{-\sigma _1-\gamma  }=x^{-1-i\zeta -\gamma }=x^{-1-\gamma }\left(\cos(\zeta \log x)+i\sin(\zeta \log x) \right)
\eqns
and we deduce, for all $t$ fixed:

\bqn
v(t, x)
&=&x^{-1-\gamma }\left(h_1(t) \cos\left(\zeta \log x\right)+h_2(t)\sin\left(\zeta \log x\right)\right)+o\left(x^{-1-\gamma } \right),\,\,\,x\to 0 \label{S6E501}\\
h_1(t)&=& \Re e A_0(t)+\Re e B_0(t) \nonumber\\
h_2(t)&=& \Im m A_0(t)-\Im m B_0(t).\nonumber
\eqn
Consider now two values of $x$:
\bqns
x_1=e^{-\frac {2 \ell \pi } {\zeta }},\,\,\,x_2=e^{-\frac {(2 \ell+1) \pi } {\zeta }}
\eqns
where $\ell\in \N$ has to be fixed. From \eqref{S6E501}:
\bqns
&&v(t, x_1)=x_1^{-1-\gamma }h_1(t)+o\left(x_1^{-1-\gamma } \right),\,\,x_1\to 0\\
&&v(t, x_2)=-x_1^{-1-\gamma }h_1(t)+o\left(x_1^{-1-\gamma } \right),\,\,x_2\to 0.
\eqns
We chose now $\ell$ large enough to have:
\bqns
&&v(t, x_1)\ge \frac {x_1^{-1-\gamma }h_1(t)} {2}>0 \,\,\,\hbox{and}\,\,\,u(t, x_2)\le- \frac {x_1^{-1-\gamma }h_1(t)} {2}<0\,\,\,\,\,\,\hbox{if}\,\,h_1(t)>0\\
&&v(t, x_1)\le \frac {x_1^{-1-\gamma }h_1(t)} {2}<0 \,\,\,\hbox{and}\,\,\,u(t, x_2)\ge - \frac {x_1^{-1-\gamma }h_1(t)} {2}>0\,\,\,\,\,\,\hbox{if}\,\,h_1(t)<0.
\eqns
\qed

\textbf{Proof of Theorem \ref{S1Th421}.} We argue by contradiction and suppose that such a local solution, that we denote  $\tilde v$, exists on some time interval $(0, T)$. We may suppose without loss of generality that $T<-\gamma ^{-1}$. By hypothesis, $\mathcal M_{\tilde v  }$ satisfies  all the assumptions in Theorem \ref{S4TU}.  By \eqref{S6E201}, $U_2$ satisfies \eqref{S4ES10} for any $\rho > 1-\gamma $. By \eqref{S6E210} and the property \eqref{S3EPO4} of $\Omega _2$,  $U_2$ also satisfies
\eqref{S4ES678} for any $\rho >1-\gamma $.
It follows by \eqref{S6E500}  that $\mathcal M _{ v }$ also  satisfies the hypothesis of Theorem \ref{S4TU} for any $\rho >1-\gamma$.  Then $\mathcal M_{\tilde v  }(t, s)$=$\mathcal M_{v}(t, s)$ for $t\in [0, T)$, $s\in \mathscr S(\rho , \rho -\gamma )$ and therefore $v(t)=\tilde v(t)$ for  $t\in [0, T)$ but this is not possible since $v$ takes  positive and negative values in $(0, \infty)$. This contradiction  concludes the proof.
\qed 

\begin{remark}
\label{S6R1}
If $\theta \in (0, 1)$ the existence of a unique global non negative solution $\mu $ is proved in Theorem 4.1 of \cite{BW}. It immediately follows from this and related results in Section 4 of \cite{BW} that the Mellin transform of $\mu $ is such that $\mathcal M _{ \mu  }\in \mathscr{C}([0, \infty); E' _{ 1+\gamma , \infty })$ and satisfies \eqref{S4ES678}-\eqref{S4ES6780} for any $\rho >1+\gamma $ and $W_0(s)=1$. We deduce by Theorem \ref{S4TU} that $\mu =v$, the solution obtained in Theorem \ref{S6T12345}, on $t\in (0, -\gamma ^{-1})$.
\end{remark}

\section{Appendix.}

\subsection{Uniqueness of bounded analytic solutions of \eqref{S2eqconst}.}
\label{Unique}
\begin{theorem}
\label{S4TU}
Given  any $T>0$ and  $W_0(s)$  a bounded and analytic function on a strip $ \mathscr{S}(\rho, \rho +|\gamma|  )$ for some $\rho >0$,  
there  exists at most one  solution $W$ to the equation \eqref{S2eqconst} for  $t\in (0, T)$, such that,  for all $t\in (0, T)$, $W(t, s)$ is analytic on the strip $ \mathscr{S}(\rho, \rho +|\gamma|  )$, satisfying
\bqn 
&&W \in C\left([0, T)\times \overline{\mathscr{S}(\rho, \rho +|\gamma|   )}\right)
\label{S4ES10}\\
&&\sup \left\{|W(t, s)|;\,\,0\le t\le T, \, s\in \overline{\mathscr{S}(\rho, \rho +|\gamma|   )} \right\}<\infty. \label{S4ES678}\\
&&W(0, s)=W_0(s),\,\,\,\forall s\in  \mathscr{S}(\rho, \rho +|\gamma|  ) \label{S4ES6780}
\eqn
\end{theorem}
\textbf{Proof.} 
Suppose that we have two solutions $W_{ \ell }(t, s)$, $\ell=1, 2$,  analytic on the strip $ \mathscr{S}(\rho, \rho +|\gamma|  )$, satisfying \eqref{S4ES10}-\eqref{S4ES6780} and  denote $W =W _1-W _2$. The function $W$  satisfies  the same conditions and $W(0)=0$. Given any $T'<T$, let $\alpha (t)$ be a $C^\infty$ cutt-off function satisfying $\alpha (t)=1$ for $0\le t\le T'$ and $\alpha (t)=0$ if $t\ge T$. If we define:
$$
\widehat W(t, s)=W(t, s)\alpha (t)
$$
we have
\bqn
\label{S7E1}
\frac {\partial \widehat W} {\partial t}(t, s)=\Phi (s)\widehat W(t, s+\gamma )+r(t, s)
\eqn
where the function $r$ is bounded in $(0, T)\times \mathscr{S}(\rho, \rho +|\gamma|  )$ and $r(t)\equiv 0$ for $0\le t\le T'$. We apply now the Laplace transform in $t$ at both sides of \eqref{S7E1} and obtain, for $\Re e(z)>0$ and $s\in \mathscr{S}(\rho, \rho +|\gamma|  )$:
\bqn
\label{S7E2}
z \widetilde W(z, s)=\Phi (s)\widetilde W(z, s+\gamma )+\widetilde r(z, s),
\eqn
where, for some constant $C>0$,
\bqn
\label{S7E3}
|\widetilde r(z, s)|\le Ce^{-T'\Re e(z)},\,\,\forall s\in\mathscr{S}(\rho, \rho +\gamma  ),\,\,\Re e(z)>0.
\eqn
By the linearity of the equation in \eqref{S7E2}  we may write $\widetilde W=\widetilde W _{ part }+\widetilde W _{ hom }$ where $\widetilde W _{ hom }$ solves
\bqn
\label{S7E4}
z \widetilde W _{ hom }(z, s)=\Phi (s)\widetilde W _{ hom }(z, s+\gamma ),\,\,\,\forall s\in\mathscr{S}(\rho, \rho +\gamma  ),\,\,\Re e(z)>0
\eqn
and $\widetilde W _{ part }$ is a particular solution of  \eqref{S7E2}. Arguing as with the function $\mathscr{U}$  defined by \eqref{S5E146} in the Proof of Proposition \eqref{S5PropU}  it follows that,  if  $\widetilde V(s)$ is the function defined by \eqref{S5E2}, then 
\bqn
\label{S7E5}
\widetilde W _{ part }(z, s)=\frac {i} {\gamma z\, \widetilde V(s)}\int _{\Re e(\sigma )=\sigma _0}\frac {(-z)^{\frac {\sigma -s} {\gamma }}\widetilde V(\sigma )\widetilde r(z, \sigma ) d\sigma } {\left(1-e^{-\frac {2i\pi } {\gamma }(s-\sigma)}\right)}
\eqn
satisfies  \eqref{S7E2} and, by \eqref{S7E3},
\bqn
\label{S7E6}
|\widetilde W _{ part }(z, s)|\le Ce^{-T'\Re e(z)},\,\,\forall s\in\mathscr{S}(\rho, \rho +|\gamma|  ),\,\,\Re e(z)>0.
\eqn
It is simplet to write our next argument if we distinguish now  the cases $\gamma >0$ and $\gamma <0$, although  the proof is completely similar in both cases. Let us then assume from now on that $\gamma >0$.  We first perform  the change of variables:
\bqn
\zeta =e^{-\frac {2i\pi } {\gamma }(s-\rho )},\,\,\,G(z, \zeta )=\widetilde W(z, s). \label{S7E6BB}
\eqn
For all $z\in \C$ such that $\Re (z)>z_0$, the function $ G(z, \zeta )$ is now analytic with respect to $\zeta $ for $\zeta \in \C\setminus \overline{\R^+}$ and bounded on $\C\setminus \overline{\R^+}$.
We also have, using that $\widetilde W \in C((0, \infty)\times \overline{ \mathscr{S}(\rho +|\gamma| )})$:
\bqn
\label{SAE432}
\begin{cases}
&\forall s=\rho +iv: \widetilde W(z, s)=G(z, x-i0), \,x:= e^{\frac {2\pi v} {\gamma }}\\
&\forall  s=\rho +\gamma +iv:  \widetilde W(z, s)=G(z, x+i0), \,x:= e^{\frac {2\pi v} {\gamma }}.
\end{cases}
\eqn
We  also define:
\bqn
\label{SAE434}
\begin{cases}
&\displaystyle{\widetilde \varphi  (\zeta )=\Phi (s)=\frac {(\rho -\sigma _1-\frac {\gamma } {2i\pi } \log \zeta )(\rho -\sigma _2-\frac {\gamma } {2i\pi }\log \zeta )}
{(\rho -\frac {\gamma } {2i\pi }\log \zeta )}},\,\,\,\forall \zeta \in \C\\
&\displaystyle{\varphi (x)=\lim _{ \varepsilon \to 0 }\widetilde \varphi  \left(xe^{-i\varepsilon } \right),\,\,\forall x>0.}
\end{cases}
\eqn
where $\log (\zeta )=\log|\zeta |+i\arg(\zeta )$, and $\arg (\zeta )\in [0, 2\pi )$.
The equation reads:
$$
G(z, x-i0)=\frac {\varphi (x)} {z}G(z, x+i0),\,\,\forall x>0.
$$
If we denote:
$$
m(z, \zeta )=\frac {1} {2i\pi }\int _0^\infty Log\left( \frac {\varphi (\lambda )} {z}\right)\left( \frac {1} {\lambda -\zeta }-\frac {1} {\lambda -\lambda _0}\right)d\lambda 
$$
where $Log (\zeta )=\log|\zeta |+i Arg(\zeta )$, and $Arg (\zeta )\in (-\pi /2, \pi/2 ]$, this is an analytic function on $\C\setminus \overline{\R^+}$ and  by Plemej-Sojoltski formulas:
$$
 \frac {\varphi (x)} {z}=\frac {e^{m(z, x+i0)}} {e^{m(z, x-i0)}},\,\,\,\forall x>0.
$$
We deduce from the equation:
$$
e^{m(z, x-i0)}G(z, x-i0)=e^{m(z, x+i0)}G(z, x+i0),\,\,\forall x>0
$$
and therefore, since  $e^{m(z, \zeta )}G(z, \zeta )$ is analytic with respect to $\zeta $ for $\zeta \in \C\setminus \R^+$,  the function 
$$
C(z, \zeta )=e^{m(z, \zeta )}G(z, \zeta )
$$
is analytic on $\C\setminus\{ 0\}$.
It remains to check the behavior of $C(z, \zeta )$ as $\zeta \to 0$ and $|\zeta |\to \infty$.

By definition:
\bqns
\varphi  (\zeta )&=&\frac {(\rho -\sigma _1-\frac {\gamma } {2i\pi }\log \zeta  )(\rho -\sigma _2-\frac {\gamma } {2i\pi }\log \zeta  )}
{(\rho -\frac {\gamma } {2i\pi }\log \zeta  )}\\
&=&-\frac {\gamma } {2i\pi }(\log |\zeta| )\frac {(1+\frac {\theta_1} {\log |\zeta| })(1+\frac {\theta_2} {\log |\zeta| })}
{(1+\frac {\theta_0} {\log |\zeta |})}
\eqns
where
$$
\theta_\ell=-\frac {2i\pi (\rho -\sigma_\ell)} {\gamma }+i \arg (\zeta ),\,\,\ell=1, 2;\,\,\,\,\theta_0=-\frac {2i\pi \rho } {\gamma }+i \arg (\zeta ).
$$
Then, as $|\zeta| \to 0$ or $|\zeta |\to \infty$
\bqn
\label{SAE786}
\varphi (\zeta )=-\frac {\gamma } {2i\pi }(\log |\zeta |)\left(1+\mathcal O\left(\frac {1} {|\log |\zeta ||}\right) \right).
\eqn
It follows that
$$
\frac {1} {2i\pi }Log\left( \frac {\varphi (\lambda )} {z}\right)=\frac {1} {2i\pi }Log\left( \frac {\gamma |\log \lambda |} {2\pi |z|}\right)
+\frac {1} {2\pi }\arg\left(-\frac {\gamma \log \lambda  } {2i\pi z} \right).
$$
Since $Arg z\in (-\pi /2, \pi /2)$ we have
$$
Arg\left(-\frac {\gamma \log \lambda  } {2i\pi z} \right)=Arg\left(i \frac {\gamma \log \lambda  } {2\pi } \right)-Arg(z)
$$
where
\[
Arg\left(i \frac {\gamma \log \lambda  } {2\pi } \right)=
\begin{dcases}
-\frac {\pi } {2}, & \hbox{if}\,\,0<\lambda <1 \\
 \frac {\pi } {2}, &  \hbox{if}\,\, \lambda >1.
\end{dcases}
\]
We may then write
\bqn
m(z, \zeta )&=&I_1(z, \zeta )+I_2(z, \zeta )+I_3(z, \zeta )\label{S6Em}\\
I_1(z, \zeta )&=&\frac {1} {2i\pi }\int _0^\infty \log\left( \frac {\gamma |\log \lambda |} {2\pi |z|}\right) \left( \frac {1} {\lambda -\zeta }-\frac {1} {\lambda -\lambda _0}\right)d\lambda \nonumber\\
I_2(z, \zeta )&=&\frac {1} {2\pi }\left(-\frac {\pi} {2}-Arg(z) \right)\int _0^1 \left( \frac {1} {\lambda -\zeta }-\frac {1} {\lambda -\lambda _0}\right)d\lambda \nonumber \\
I_3(z, \zeta )&=&\frac {1} {2\pi }\left(\frac {\pi} {2}-Arg(z) \right)\int _1^\infty \left( \frac {1} {\lambda -\zeta }-\frac {1} {\lambda -\lambda _0}\right)d\lambda \nonumber
\eqn
If we take $\lambda _0=i$:
\bqns
\int _0^1 \left( \frac {1} {\lambda -\zeta }-\frac {1} {\lambda -\lambda _0}\right)d\lambda=-\frac {i\pi } {4}-\frac {Log 2} {2}+Log\left(\frac {\zeta -1} {\zeta } \right)\\
\int _1^\infty \left( \frac {1} {\lambda -\zeta }-\frac {1} {\lambda -\lambda _0}\right)d\lambda=-\frac {i\pi } {4}+\frac {Log 4} {4}-Log\left(1-\zeta  \right)
\eqns
and then,
\bqns
(I_2+I_3)(z, \zeta )&=&\frac {1} {2\pi }\left\{ \left( \frac {\pi } {2}+Arg (z)\right)Log\left(\frac {\zeta } {\zeta-1 } \right)-\left( \frac {\pi } {2}-Arg (z)\right)Log\left(1-\zeta \right)\right\}+R(z)\\
R(z)&=&\frac {1} {2\pi }\left\{ -\left( \frac {\pi } {2}+Arg (z)\right)\left(\frac {i\pi } {4}+\frac {Log 2} {2} \right)+\left( \frac {\pi } {2}-Arg (z)\right)\left( \frac {Log 2} {2}-\frac {i\pi } {4}\right)\right\}.
\eqns
As $|\zeta |\to 0$, 
\bqns
Log\left(\frac {\zeta } {\zeta-1 } \right)&=&Log(\zeta )-\log(\zeta -1)=log|\zeta |+iArg(\zeta )-Log(\zeta -1)\\
&=&\log|\zeta |\left(1+\frac {iArg(\zeta )-\log(\zeta -1)} {log|\zeta |} \right)=\log|\zeta |\left(1+\mathcal O\left( \frac {1} {\log|\zeta |}\right)\right),  |\zeta |\to 0\\
&&\Longrightarrow (I_2+I_3)(z, \zeta )= \frac {1} {2\pi }\left( \frac {\pi } {2}+Arg (z)\right)\left(\log|\zeta |+\mathcal O\left( \frac {1} {\log|\zeta |}\right)\right),\,\, |\zeta |\to 0.
\eqns

As $ |\zeta |\to \infty$
\bqns
Log (1-\zeta )&=&\log|\zeta -1|+iArg (1-\zeta )=\log|\zeta |+\log\left|1-\frac {1} {\zeta } \right|+iArg (1-\zeta )\\
&=&\log|\zeta |\left(1+\frac {\log\left|1-\frac {1} {\zeta } \right|+iArg (1-\zeta )} {\log|\zeta |} \right)=\log|\zeta |\left(1+\mathcal O\left( \frac {1} {\log|\zeta |}\right)\right)\\
&&\Longrightarrow (I_2+I_3)(z, \zeta )= -\frac {1} {2\pi }\left( \frac {\pi } {2}-Arg (z)\right)\left(\log|\zeta |+\mathcal O\left( \frac {1} {\log|\zeta |}\right)\right),\,\, |\zeta |\to \infty.
\eqns
By \eqref{S6Em}, $e^{m(z, \zeta )}=e^{I_1(z, \zeta )}e^{(I_2+I_3)(z, \zeta )}$.  We notice that $\left|e^{I_1(z, \zeta )}\right|=1$ and, as $|\zeta |\to 0$:
\bqn
\left|e^{(I_2+I_3)(z, \zeta )}\right|&=&e^{ \frac {1} {2\pi }\left( \frac {\pi } {2}+Arg (z)\right)\left(\log|\zeta |+\mathcal O\left( \frac {1} {\log|\zeta |}\right)\right)}
\le e^{ \frac {1} {2\pi }\left( \frac {\pi } {2}+Arg (z)\right)\left(\log|\zeta |+1\right)}\nonumber\\
&=&e^{ \frac {1} {2\pi }\left( \frac {\pi } {2}+Arg (z)\right)}
e^{ \frac {1} {2\pi }\left( \frac {\pi } {2}+Arg (z)\right)\log|\zeta |}=e^{ \frac {1} {2\pi }\left( \frac {\pi } {2}+Arg (z)\right)} |\zeta |^{ \frac {1} {2\pi }\left( \frac {\pi } {2}+Arg (z)\right)} \label{SAEarg1}
\eqn
On the other hand, as $|\zeta |\to \infty$, a similar argument gives:
\bqn
\label{SAEarg2}
\left|e^{(I_2+I_3)(z, \zeta )}\right|\le e^{ -\frac {1} {2\pi }\left( \frac {\pi } {2}-Arg (z)\right)} |\zeta |^{- \frac {1} {2\pi }\left( \frac {\pi } {2}-Arg (z)\right)}.
\eqn
Since $Arg(z)\in (-\pi /2, \pi /2)$:
\bqn
\label{SAEarg}
0<\frac {1} {2\pi }\left( \frac {\pi } {2}+Arg (z)\right)<\frac {1} {2}.
\eqn
From the boundedness of the function $G(z,\cdot)$  on $\C$ and  \eqref{SAEarg1}-\eqref{SAEarg}, we deduce  that for all $z\in\C, \, \Re e(z)>z_0$, the function
$C(z, \zeta )=G(z, \zeta )e^{m(z, \zeta )}$  is bounded as $|\zeta |\to 0$ and $|\zeta |\to \infty$. It follows that $C(z, \cdot)$ is independent of $\zeta $. Using \eqref {SAEarg1} again
$$
\lim _{ \zeta \to 0 } C(z, \zeta )=0
$$
and we deduce that $C(z, \zeta )=0$ for all $\zeta $, then $G(z)\equiv 0$ for all $z\in\C, \, \Re e(z)>z_0$. Therefore $\widetilde W _{ hom }=0$  and then $\widetilde W=
\widetilde W _{ part }$. Laplace's inversion then yields:
\bqn
\label{SAE421}
\widehat W(t, s)=\frac {1} {2i\pi }\int_{b-i\infty}^{b+i \infty} \widehat W _{ part }(t, s)e^{z t}dz
\eqn
for any $b>0$. Then, \eqref{S7E6} implies  $\widehat W(t, s)=W(t, s)=0$ for all $0\le t\le T'$ and $s\in \mathscr{S}(\rho, \rho +|\gamma|  )$.
\qed
\vskip0.15cm 
It is not always possible to apply Theorem \ref{S4TU} to the solutions  of a Cauchy problem associated to \eqref{S2eqconst}. That is the case when $\gamma >0$ and  consider the solution $U$, obtained in Section \ref{Extension}, Proposition \ref{S5PropU}.  Our next result is then useful:
\begin{theorem}
\label{S4TU2}
Suppose  $\gamma >0$. Given  any $T>0$ and  $W_0(s)$  such that $sW_0(s)$ is a bounded and analytic function on the strip $ \mathscr{S}(-\gamma , \varepsilon  )$ for some $\varepsilon  >0$,  
there  exists at most one  solution $W$ to the equation \eqref{S2eqconst} for  $t\in (0, T)$, such that,  for all $t\in (0, T)$, $sW(t, s)$ is analytic on the strip $ \mathscr{S}(-\gamma , \varepsilon  )$, satisfying
\bqns 
&&sW \in C\left([0, T)\times \mathscr{S}(-\gamma +\delta , \varepsilon -\delta  )\right),\,\,\hbox{for some}\,\delta \in (0, \varepsilon )
\label{S4ES1067}\\
&&\sup \left\{|sW(t, s)|;\,\,0\le t\le T, \, s\in \mathscr{S}(-\gamma +\delta , \varepsilon -\delta ) \right\}<\infty. \label{S4ES67867}\\
&&W(0, s)=W_0(s),\,\,\,\forall s\in  \mathscr{S}(-\gamma , \varepsilon ) \label{S4ES678067}
\eqns
\end{theorem}
\textbf{Proof of Theorem \ref{S4TU2}.} Assume  the existence of two such solutions to \eqref{S2eqconst} and call $W$ their difference. Then, we  define the  function:
\bqns
\label{SAEH}
H(t, s)=sW(t, s)
\eqns
If $\rho >0$ is such that $(\rho , \rho +\gamma )\subset (-\gamma , \varepsilon )$, by our hypothesis on $W$:
\bqns
\frac {\partial H(t, s)} {\partial t}=\frac {(s-\sigma _1)(s-\sigma _2)} {s+\gamma }H(t, s+\gamma)\\
H(0, s)=sW_0(s),\,\,\forall s \in \mathscr{S}(\rho, \rho +\gamma  ).
\eqns
The proof follows now  the same arguments used in the proof of Theorem \ref{S4TU}, applying  to $H$ the same arguments used  in the proof of Theorem \ref{S4TU} with $W$. This amounts just to consider the new function 
$$
\Psi (s)=\frac {(s-\sigma _1)(s-\sigma _2)} {s+\gamma }
$$
instead of $\Phi $. 

We first consider the case where $T<\gamma ^{-1}$. Then,  all  the beginning of the proof of Theorem \ref{S4TU} may be exactly reproduced  until the formula \eqref{S7E4}, with the interval $(\rho , \rho +\gamma )$. In order  to obtain a particular solution $\widetilde H _{ part }$  of 
\bqn
\label{S7E2HH}
z \widetilde W(z, s)=\Psi (s)\widetilde W(z, s+\gamma )+\widetilde r(z, s),
\eqn
we consider the function:
\bqns
\mathscr{V}(s)=\frac {(-\gamma) ^{\frac {s} {\gamma }}\Gamma \left(\frac {s-\sigma _1} {\gamma } \right)\Gamma \left(\frac {s-\sigma _2} {\gamma } \right)} {\Gamma \left(1+\frac {s} {\gamma } \right)}.
\eqns
It is  straightforward to check that $\mathscr{V}(s)$ satisfies:
\bqns
\mathscr{V}(s+\gamma )=-\Psi (s)\mathscr{V}(s),\,\,\,s\in \C\setminus \{s\in \C; s=\sigma _1-m\gamma ,\, m=0, 1, 2, \cdots\}
\eqns
and, arguing as for the function $V$ in \eqref{S5EV}, its behavior as $|\Im m (s)|\to \infty$ with $\Re e(s)$ bounded is  such that the function:
$$
\widetilde H _{ part }(t, s)=\frac {i} {\gamma z\, \mathscr{V}(s)}\int _{\Re e(\sigma )=\sigma _0}\frac {(-z)^{\frac {\sigma -s} {\gamma }}\mathscr{V}(\sigma )\widetilde r(z, \sigma ) d\sigma } {\left(1-e^{-\frac {2i\pi } {\gamma }(s-\sigma)}\right)}
$$
satisfies the equation  \eqref{S7E2HH}  and the estimate \eqref{S7E6} for $t\in (0, T)$.

The argument for $\widetilde H_{ hom }$ is now very similar using the new functions $\mathscr{G}(t, \zeta )=H(t, s)$ instead of $G$ in \eqref{S7E6BB}, and $\tilde \psi (\zeta )$, $\psi (x)$ instead of $\tilde \varphi (\zeta )$,  $\varphi  (x)$ in \eqref{SAE432}, \eqref{SAE434}. Since $\psi (\zeta )$ may still be estimated as in \eqref{SAE786} the end of the argument follows straightforwardly in the same way to prove that $H(t)\equiv 0$ for $t\in [0, T)$. This proves Theorem \ref{S4TU2} if $T<\gamma ^{-1}$.

The result for $T>\gamma ^{-1}$ follows by iteration of the uniqueness of solutions on intervals $[m\gamma /2, (m+1)\gamma /2)$ for $m =0, 1, 2,\cdots, M$ where
$M=[2T/\gamma] $, the integer part of $2T/\gamma $, and finally on $[M\gamma /2, T)$.
\qed
\section*{Acknowledgements}
The research of the author is supported by  grants  MTM2014-52347-C2-1-R of DGES and  IT641-13 of the Basque Government. The hospitality of IAM at the University of Bonn, where part of this work was done,  is  warmly acknowledged. The author is grateful to Pr. J. Bertoin for enlightening comments on the results in \cite{BW}.

 \end{document}